\documentclass[ASNA,twocolumn]{USG}

\usepackage[capitalize,nameinlink]{cleveref}
\crefname{figure}{Fig.}{Figs.}
\Crefname{figure}{Fig.}{Figs.}
\crefname{table}{Table}{Tables}
\Crefname{table}{Table}{Tables}

\graphicspath{{./}{./images/}{./experiments/}{./figures/}{./figures/plots/}}

\articletype{Research Article}
\journal{International Journal for Numerical Methods in Engineering}
\received{}
\revised{}
\accepted{}
\volume{0}
\copyyear{2026}
\startpage{1}
\articledoi{10.1002/nme.0000}

\begin{document}


\title{FEVessel: Mesh-Independent Analysis of 3D Pressure Vessels with the Label-Free Pretrained Finite Element Method}

\author[1]{Yipin Sun}[https://orcid.org/0009-0003-7752-8321]
\author[2]{Yizheng Wang}[https://orcid.org/0000-0002-3899-7008]
\author[2]{Yuzhou Lin}[https://orcid.org/0009-0006-5541-3875]
\author[1]{Baiyang Zheng}[https://orcid.org/0000-0001-8616-6180]
\author[4]{Xiaoying Zhuang}[https://orcid.org/0000-0001-6562-2618]
\author[1,3]{Timon Rabczuk}[https://orcid.org/0000-0002-7150-296X]

\authormark{Sun \textsc{et al.}}
\titlemark{FEVessel: Mesh-Independent Analysis of 3D Pressure Vessels}

\address[1]{\orgdiv{Institute of Computational Mechanics x AI}, \orgname{Fudan University},
\orgaddress{\state{Shanghai 200090}, \country{China}}}

\address[2]{\orgdiv{Department of Engineering Mechanics}, \orgname{Tsinghua University},
\orgaddress{\state{Beijing 100084}, \country{China}}}

\address[3]{\orgdiv{Institute of Structural Mechanics}, \orgname{Bauhaus University Weimar},
\orgaddress{\state{99423 Weimar}, \country{Germany}}}

\address[4]{\orgdiv{Institute of Photonics, Faculty of Mathematics and Physics}, \orgname{Leibniz University Hannover},
\orgaddress{\state{30167 Hannover}, \country{Germany}}}

\corres{%
Yizheng Wang, \email{wang-yz19@tsinghua.org.cn}\\
Xiaoying Zhuang, \email{zhuang@iop.uni-hannover.de}\\
Timon Rabczuk, \email{timonrabczuk@fudan.edu.cn}}

\fundingInfo{China National Scholarship}


\abstract[Abstract]{%
Pressure vessel analysis in the chemical, nuclear, and
new-energy industries requires solving the same elasticity
problem across many materials, geometries, and loads, where
mesh quality and repeated solving govern both accuracy and cost.
The finite element method (FEM) cannot amortise this repeated
cost and fails on degenerate meshes, while the neural operators
meant to replace it still need labelled data that FEM must
generate.
This paper proposes FEVessel, an adaptation of the Pretrained
Finite Element Method (PFEM) to three-dimensional (3D) pressure
vessels, and validates four capabilities across the two
limitations above.
FEVessel i) encodes each vessel as a point cloud with
coordinate, material, and load channels, ii) pretrains a
Transolver operator on the total potential energy instead of FEM
labels, and iii) warm-starts iterative solvers with its
prediction.
A single model generalises across material, geometry, and
boundary conditions at a $1.35\%$ relative displacement error,
and its $2.07\%$ strain error is about $4.7$ times lower than
that of a supervised Fourier neural operator ($9.72\%$), whose
structured grid cannot preserve the through-thickness strain.
Its warm start cuts algebraic multigrid iterations from $195$ to
$18$, a $9.2\times$ end-to-end wall-clock speedup at the
$10^{-3}$ engineering tolerance.
The model transfers across mesh resolutions without retraining,
holding about $3\%$ error at only $30\%$ of the training point
density.
On inverted and sliver meshes where FEM fails, the error remains
below $3.66\%$.
To our knowledge, this is the first systematic study of
mesh-independent solution on industrially relevant 3D pressure
vessels with degenerate meshes.
Because training needs no labels, FEVessel works exactly where
FEM cannot supply any, removing manual mesh repair from the
analysis pipeline.}

\maketitle

\section{Introduction}\label{sec:intro}

Pressure vessels are widely used for high pressure fluid
storage in the chemical, petroleum, nuclear, and new-energy
industries~\cite{moss2004pressure,peters2023overview}.
Due to the extremely high internal pressure, the failure of a pressure
vessel can cause catastrophic consequences, including loss of
life and significant economic damage~\cite{alves2025probabilistic,li2025reliability}.
Therefore, the structural integrity of the pressure vessel is critical throughout its service life. 
To ensure this integrity, the governing linear elasticity equations must be solved accurately
across the parameter space defined by material, geometry, and loading conditions~\cite{johnson2023parametric,zhang2022surrogate,kang2025condition,wang2024dcem}.

The linear elasticity of a pressure vessel is governed by partial differential equations (PDEs). 
The finite element method (FEM) has long been the standard tool for solving these PDEs, owing to its
rigorous mathematical foundation and controllable accuracy~\cite{zienkiewicz1977finite,hughes2003finite}.
For instance, FEM has been applied to resolve stress
concentrations at nozzle junctions in pressure vessels, where
conventional analytical formulae are no longer applicable~\cite{bozkurt2021comparison}.
However, FEM solves each configuration independently: any change in the material, geometry, or boundary conditions requires solving the large stiffness system again,
incurring a high computational cost~\cite{wang2026pretrain}.
This cost has driven scholars to develop dedicated remedies. 
Kudela et al.~\cite{kudela2022recent} surveyed surrogate models that
amortize repeated FEM evaluations across parametric studies and optimization, 
while Sun et al.~\cite{sun2025novel} proposed an integrated surrogate model to  replace the costly finite-element evaluations. 
Although these studies improve the efficiency of repeated solving, 
FEM still relies heavily on mesh resolution and element quality in the linear elasticity 
analysis of pressure vessels.
Specifically, when the mesh resolution is sufficiently refined, FEM produces reliable results; 
however, if the resolution becomes insufficient, the computed results are no longer reliable.
For example, Flanagan et al.~\cite{flanagan2025through} found that a mesh with too few
elements through the thickness cannot capture the stress variation, leaving the result unreliable. 
Even a sufficiently refined mesh can fail on element quality. The automated 
meshing of computer-aided design (CAD) geometry still produces inverted or sliver elements, which defeat the
FEM solver~\cite{satheesh2024novel,van2022incorporation,wei2024structural}. 
Sukumar et al.~\cite{sukumar2022virtual} noted that even robust tetrahedral
meshers cannot remove these degenerate elements. Because an inverted element
has a negative Jacobian determinant, the shape-function gradients
cannot be evaluated and the assembly fails before any system is
formed~\cite{prabhune2024isoparametric}. A sliver element, by
contrast, keeps a positive Jacobian, so the assembly succeeds, but the
stiffness matrix is ill-conditioned and the solver fails to
converge~\cite{conley2016overcoming}. Either way FEM produces no
solution and the analyst is forced into a costly manual remeshing
step. Both limitations stem from FEM's reliance on a high-quality, well-resolved mesh,
and therefore motivate methods whose accuracy does not depend on mesh resolution or element quality.

A different paradigm solves these PDEs with neural networks rather
than the FEM, broadly termed artificial intelligence for PDEs (AI
for PDEs)~\cite{karniadakis2021physics,cuomo2022scientific}. Work in this area falls into 
three families that differ in how the network learns the solution.
(i)~The physics-informed neural network (PINN) solves these PDEs
without any FEM reference solutions. The reason is that the PINN
takes the PDE residual itself as its training loss, which lets
the governing physics constrain the network in place of
labelled data. For instance, Raissi et al.~\cite{raissi2019physics}
introduced the PINN and demonstrated it on both forward and
inverse problems governed by nonlinear PDEs across fluid,
quantum, and reaction-diffusion systems. However, each PINN
must be retrained whenever the material, geometry, or boundary
changes, so like the FEM it cannot amortise the parameter sweep. 
(ii)~Operator learning instead predicts an entire problem class
in a single forward pass and therefore amortises the parameter
sweep that the PINN cannot. The reason is that operator learning
targets a mapping between function spaces rather than a single
solution field, so one trained operator spans the entire class
without retraining. For instance, Lu et al.~\cite{lu2021learning}
established this paradigm with DeepONet, and Li
et al.~\cite{li2021fourier} then introduced the Fourier neural operator (FNO), which runs up to three orders of magnitude
faster than conventional solvers. However, most operator learning methods are
trained on FEM reference solutions, and generating these
reference solutions across the parameter space is itself
expensive, so where the FEM yields no reference solution the
operator cannot be trained. 
(iii)~Physics-informed neural operators (PINO) predict in a single
forward pass like other neural operators while embedding physical
constraints into the neural operators objective, removing the dependence on FEM reference solutions. 
For instance, Li et al.~\cite{li2024physics} combined training data with PDE residuals
in the PINO, and Eshaghi et al.~\cite{eshaghi2025variational} replaced the residual with an
energy functional in the Variational Physics-Informed Neural Operator (VINO). 
The Pretrained Finite Element Method (PFEM) of Wang et al.~\cite{wang2026pretrain}
follows this energy-functional formulation and is highly
effective for elasticity problems. PFEM uses the
Transolver backbone of Wu et al.~\cite{Wu2024TransolverAF} and
trains on unstructured point clouds, requiring no FEM reference solutions. 
PFEM also predicts the displacement field rapidly in a single forward pass. The predicted
field can in turn serve as the warm-start initial guess for a
classical iterative solver, substantially reducing the computational
cost. However, PFEM has not yet been applied to three-dimensional (3D) pressure vessels.

We close this gap with FEVessel, an adaptation of PFEM to 3D pressure vessel analysis. 
Each vessel is represented as an unstructured point cloud, and the
operator is trained from a physics-based energy functional, with no
FEM reference solutions. We validate four capabilities, the first
two addressing repeated solving across configurations and the last two mesh
dependence.

\begin{itemize}
  \item \textbf{Simultaneous M/G/B generalisation (\cref{subsec:mgb}).}
    A single pretrained model predicts the displacement field across
    material (M), geometry (G), and boundary (B) at once, reaching a $1.35\%$
    relative $L^{2}$ error without any labels, and a $2.07\%$ strain
    ($H^{1}$) error about $4.7\times$ lower than that of a supervised
    FNO trained on the same data.
  \item \textbf{Warm-start acceleration (\cref{subsec:warmstart}).} The
    same prediction, used as the initial guess for a classical
    iterative solver, cuts the algebraic multigrid iteration count from $195$ to $18$
    at an engineering tolerance of $10^{-3}$, an end-to-end speedup of
    $9.2\times$.
  \item \textbf{Cross-resolution prediction
    (\cref{subsec:resolution}).} A model trained at one resolution
    predicts directly at others without retraining, holding about
    $3\%$ mean relative $L^{2}$ error even at $30\%$ of the training
    point density.
  \item \textbf{Solvability on degenerate meshes
    (\cref{subsec:degenerate}).} On meshes that defeat a conventional
    FEM solver, FEVessel still returns physically consistent displacement
    fields across $1\%$ to $75\%$ degeneration, unaffected by inverted
    elements and holding the sliver-case mean error to $3.66\%$ at
    worst.
\end{itemize}

The remainder of the paper is organised as follows.
\Cref{sec:prerequisite} reviews the PFEM framework. \Cref{sec:method}
describes FEVessel, the adaptation for pressure vessel analysis.
\Cref{sec:experiments} presents the four sets of experiments.
\Cref{sec:discussion} discusses the engineering implications and
limitations. \Cref{sec:conclusion} concludes.

\section{Prerequisite: the PFEM framework}\label{sec:prerequisite}

%
%

We briefly review PFEM~\cite{wang2026pretrain}, the methodological
backbone of this work, focusing on physics-informed neural operators
(\cref{subsec:operator-learning}), geometric input representations
(\cref{subsec:geometric-representations}), the Transolver backbone
(\cref{subsec:transolver}), and the PFEM two-stage framework
(\cref{subsec:pfem}).


\subsection{Physics-informed neural operators}
\label{subsec:operator-learning}

Neural operators serve as fast surrogate solvers for parametric
partial differential equations. Because a neural operator
learns a mapping between function spaces, a single trained operator
predicts the solution for a new configuration directly and avoids
solving the problem again. For a parametric boundary-value problem
(BVP), the solution map is written as
\begin{equation}\label{eq:operator}
  \mathcal{H}: \mathcal{A} \to \mathcal{U}, \qquad
  a \mapsto u,
\end{equation}
where $\mathcal{A}$ and $\mathcal{U}$ are the configuration and
solution function spaces. The configuration $a \in \mathcal{A}$
encodes the geometry, material fields, and boundary data, and the
operator $\mathcal{H}$ maps the configuration to the corresponding
solution field $u \in \mathcal{U}$. The configuration $a$ and the
solution $u$ are sampled at $N$ points, and the map is approximated
by a neural operator $\mathcal{H}_\theta$ with learnable parameters
$\theta$.

Among existing architectures, DeepONet~\cite{lu2021learning} and
FNO~\cite{li2021fourier} are two representative data-driven neural
operators. DeepONet encodes the input function at a fixed set of
sensor points with a branch network, encodes the query coordinate
with a trunk network, and forms the output from the inner product
of the two networks. FNO instead defines the functions on a
structured grid and parameterises the integral kernel in the
frequency domain through the fast Fourier transform. Both methods
are trained on precomputed high-fidelity reference solutions, and
generating these solutions across the full parameter space is
expensive, so where the FEM yields no reference solution the
operator cannot be trained either.

To remove this dependence on labelled data, the PINO
embeds the governing equations into operator training. Like a
data-driven operator, PINO learns the function-space map of
\cref{eq:operator} and predicts a new configuration in a single
forward pass. Its training objective, however, is the residual
of the governing equations at collocation points, rather than a
fit to reference displacements. This residual is computed from
the predicted field and the equations alone, so it needs no
knowledge of the true solution.


\subsection{Geometric input representations}
\label{subsec:geometric-representations}

\begin{figure*}[!t]
  \centering
  \includegraphics[width=\linewidth]%
    {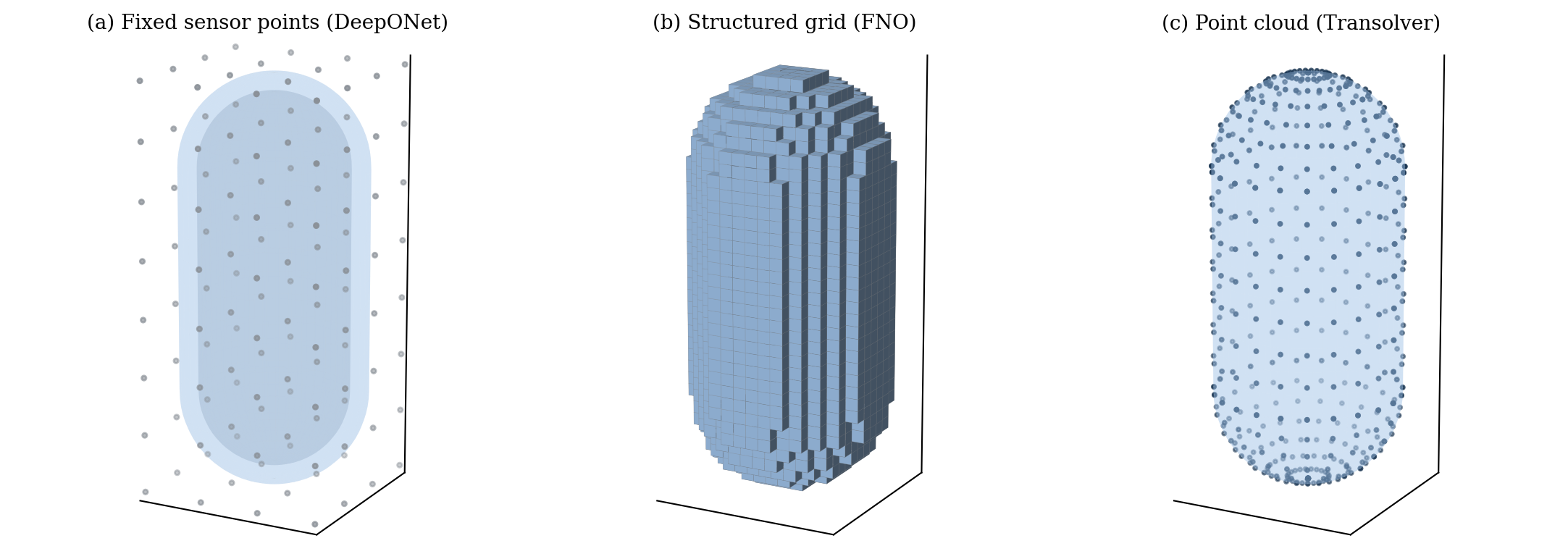}
  \caption{Geometric input representations of three neural
    operator architectures applied to a pressure vessel.
    (a)~DeepONet discretises the input function at a fixed set
    of sensor points. (b)~FNO maps the domain onto a structured
    grid (shown: voxelisation). (c)~Transolver represents the
    domain as an unstructured point cloud.}
  \label{fig:representation}
\end{figure*}

Beyond the training signal, the geometric representation of the
input constitutes a further independent limitation.
\Cref{fig:representation} contrasts the three representations
adopted by DeepONet, FNO, and Transolver. DeepONet encodes the
input function through a branch network evaluated at a fixed set
of sensor points (\cref{fig:representation}(a)), enabling
solution prediction at arbitrary query coordinates without
binding the output to a discrete mesh. For a pressure vessel,
however, the sensor layout remains fixed across the entire
parameter space, and when the inner radius $R$ or the
wall-thickness ratio $t/R$ varies the layout does not adapt,
so the input representation cannot cover the joint material,
geometry, and boundary parameter space.

FNO parameterises the integral kernel by mapping the input
onto a regular grid and applying the fast Fourier
transform (\cref{fig:representation}(b)), enabling global
feature capture in a single forward pass. The structured grid,
however, never coincides with the analysis mesh. A naive
voxelisation approximates the curved inner and outer walls as
staircase boundaries, and even a body-fitted grid must exchange
every field with the mesh by interpolation, a transfer that
degrades precisely the through-thickness displacement gradient
where the stress reaches its peak.

A point cloud places the sample points directly on the true
geometric surfaces (\cref{fig:representation}(c)), so the
curved walls and the internal-pressure boundary are represented
exactly without geometric approximation. Pressure vessel
analysis therefore calls for an operator that acts directly
on point clouds, and the Transolver architecture adopted in
this work is presented in \cref{subsec:transolver}.

\subsection{Transolver backbone}\label{subsec:transolver}

The point-cloud operator adopted by PFEM is
Transolver~\cite{Wu2024TransolverAF}, a transformer-based
architecture. Applying self-attention to every point would cost
$\mathcal{O}(N^{2})$. Transolver instead applies attention to a
small set of $S$ physics-aware tokens and reduces the cost to
$\mathcal{O}(N + S^{2})$ with $S \ll N$, so Transolver stays
efficient on the large point clouds produced by 3D meshes.

The core of Transolver is the physics-attention module. Each point
carries a tuple $\{\boldsymbol{x}^{(i)}, m^{(i)}, b^{(i)}\}$ of
coordinate, material, and boundary-condition channels, which a multi-layer perceptron (MLP)
with layer normalisation maps to a pointwise feature representation
$\mathbf{X} \in \mathbb{R}^{N \times C}$ with $C$ feature channels.
Two MLPs then map $\mathbf{X}$ to the value features
$\mathbf{U} \in \mathbb{R}^{N \times C}$ and to the slice weights
$\mathbf{M} \in \mathbb{R}^{N \times S}$, where $\mathbf{M}$ is
softmax-normalised along the token dimension, so that each row of
$\mathbf{M}$ gives a point's membership probabilities over the $S$
tokens. The physics-aware tokens
$\mathbf{Z} \in \mathbb{R}^{S \times C}$ are then obtained by
weighted aggregation,
\begin{equation}\label{eq:slice}
  Z_{JK} =
  \frac{\sum_{I=1}^{N} M_{IJ}\, U_{IK}}
       {\sum_{I=1}^{N} M_{IJ}},
\end{equation}
where $I = 1, \dots, N$ indexes the points, $J = 1, \dots, S$ the
tokens, and $K = 1, \dots, C$ the channels; that is, each token is
the membership-weighted average of the point features. Since
$S \ll N$, attention on $\mathbf{Z}$ can then be computed at low
cost. Self-attention is then applied along the token dimension. Here
$\mathbf{Q}, \mathbf{K}, \mathbf{V} \in \mathbb{R}^{S \times C}$ are
the query, key, and value matrices produced from $\mathbf{Z}$, and
the tokens are updated as
\begin{equation}\label{eq:tokenattn}
  \mathbf{Q}, \mathbf{K}, \mathbf{V} = \mathrm{MLP}(\mathbf{Z}),
  \qquad
  \mathbf{Z}^{t}
  = \mathrm{softmax}\bigl(\mathbf{Q}\mathbf{K}^{\top}\bigr)
    \mathbf{V},
\end{equation}
giving the transformed tokens $\mathbf{Z}^{t} \in \mathbb{R}^{S
\times C}$, in which long-range correlations between the physical
regions are captured. The pointwise output is recovered by
deslicing with the same weights $\mathbf{M}$,
\begin{equation}\label{eq:deslice}
  X^{t}_{IK} = \sum_{J=1}^{S} M_{IJ}\, Z^{t}_{JK},
\end{equation}
so that each point inherits the updated token states in proportion
to its region membership, producing the layer output $\mathbf{X}^{t}
\in \mathbb{R}^{N \times C}$ in the same pointwise format as the
input. Multiple such Transolver layers can be stacked to improve
expressiveness.

Crucially, the Transolver backbone operates entirely on point clouds
and does not depend on element connectivity, which underlies the
mesh-independent prediction capability evaluated in
\cref{subsec:degenerate}.

\subsection{The PFEM two-stage framework}\label{subsec:pfem}

PFEM consists of two stages: a physics-informed pretraining
stage and an optional warm-start stage.
Both stages rely exclusively on the governing physical equations and
do not require labelled solution data.

\subsubsection{Pretraining stage}\label{subsubsec:pretraining}

In the pretraining stage, the Transolver-based neural operator
$\mathcal{H}_\theta$ is trained from the governing equations alone.
Consider a parametric BVP
\begin{equation}\label{eq:bvp}
  \mathcal{P}[u(\boldsymbol{x}); a]
  = f(\boldsymbol{x}; a) \;\text{in } \Omega(a),
  \qquad
  \mathcal{B}[u(\boldsymbol{x}); a]
  = g(\boldsymbol{x}; a) \;\text{on } \partial\Omega(a),
\end{equation}
where $\boldsymbol{x}$ is the spatial coordinate, $\Omega(a)$ is the
configuration-dependent problem domain with boundary
$\partial\Omega(a)$, $\mathcal{P}$ and $\mathcal{B}$ are the domain and
boundary operators, and $f$, $g$ are the prescribed source and boundary
data.
A strong-form (PINO-style~\cite{li2024physics}) loss, which
penalises the pointwise residuals of \cref{eq:bvp} at collocation
points, applies to general, possibly non-self-adjoint PDEs. The
present problem (linear elasticity) is self-adjoint and admits a
total-potential-energy functional $\Pi[u; a]$ whose Euler--Lagrange
equation is \cref{eq:bvp}, so we adopt the energy-form
(VINO-style~\cite{eshaghi2025variational}) loss
\begin{equation}\label{eq:pfem-loss-energy}
  \mathcal{L}_{\text{PFEM}}^{\text{e}}(\theta)
  = \frac{1}{|\mathcal{D}|}
    \sum_{a \in \mathcal{D}}
    \Pi\bigl[\mathcal{H}_\theta(a);\, a\bigr],
\end{equation}
which trains the network to minimise the total potential energy over
the training set $\mathcal{D}$. The explicit $\Pi$ for the
pressure-vessel problem is given in \cref{subsec:problem}.

The energy form is preferred for two reasons. It uses derivatives of
$u$ one order lower than the strong form (first-order strains, not
second-order divergences), improving training stability; and it
embeds the natural (Neumann) boundary condition directly in the
boundary integral of $\Pi$, so no separate boundary-loss term is
needed. The essential (Dirichlet) condition is not enforced by $\Pi$
and is imposed separately by hard-constraining the network output on
the Dirichlet boundary $\Gamma_u$ (\cref{subsec:encoding}).

\subsubsection{Warm-start stage}\label{subsubsec:warmstart}

In the warm-start stage, the pretrained prediction is used as the
initial guess for a classical iterative solver applied to the
assembled stiffness system $\mathbf{K}\mathbf{U} = \mathbf{F}$.
Writing $\mathbf{U}^{\mathrm{NO}} \equiv \mathcal{H}_\theta(a)$ for
the neural-operator prediction at the mesh nodes, the iteration reads
\begin{equation}\label{eq:warmstart}
  \mathbf{U}^{(k+1)} = \boldsymbol{\phi}\bigl(
    \mathbf{U}^{(k)};\, \mathbf{K},\, \mathbf{F}
  \bigr),
  \qquad
  \mathbf{U}^{(0)} = \mathbf{U}^{\mathrm{NO}},
\end{equation}
where $\mathbf{U}^{(k)}$ is the nodal displacement vector at
iteration $k$, $\mathbf{K}$ and $\mathbf{F}$ are the global stiffness
matrix and nodal load vector, and $\boldsymbol{\phi}(\cdot)$ is the
chosen iterative algorithm (e.g., conjugate gradient or algebraic
multigrid). The iteration terminates once the relative residual of
the unbalanced nodal force $\mathbf{r}(\mathbf{U}) =
\mathbf{K}\mathbf{U} - \mathbf{F}$ drops below the tolerance,
\begin{equation}\label{eq:residual}
  \frac{\bigl\lVert \mathbf{r}\bigl(\mathbf{U}^{(k+1)}\bigr)
    \bigr\rVert}{\lVert \mathbf{F} \rVert} < \mathrm{tol},
\end{equation}
where $\lVert\cdot\rVert$ is the Euclidean norm and $\mathrm{tol}$ is
the prescribed tolerance. Because $\mathbf{U}^{\mathrm{NO}}$ is
already physics-consistent and close to the converged solution, far
fewer iterations are required than from a zero or random initial
guess. \Cref{app:warmstart-theory} makes this precise for the linear
pressure vessel system. The iteration count derived in
\cref{eq:app-itercount} grows only with the logarithm of the
initial relative residual. A warm start that lowers that residual
by a larger factor therefore removes more iterations. This saving
is fixed and tolerance-independent, given in closed form by
\cref{eq:app-deltan}.

\section{FEVessel: PFEM for pressure vessel
analysis}\label{sec:method}

\begin{figure*}[!t]
  \centering
  \includegraphics[width=\linewidth]{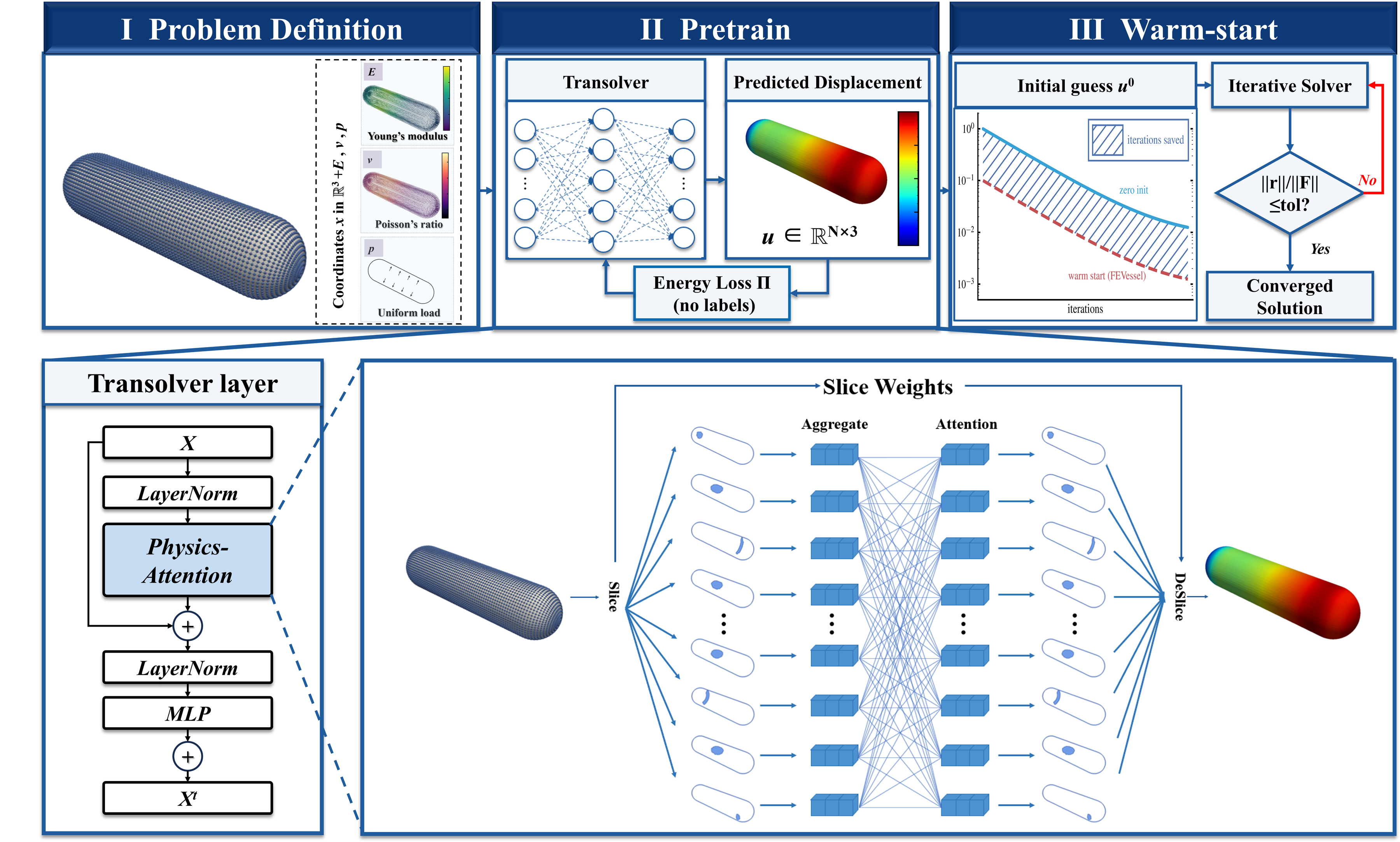}
  \caption{Overview of the FEVessel input and output for
    3D pressure vessel analysis. Stage~I defines the
    problem. Each mesh node carries its coordinates
    $\boldsymbol{x}\in\mathbb{R}^{3}$ together with $E$, $\nu$, and $p$.
    Stage~II pretrains the Transolver operator, which maps this point
    cloud to the displacement field
    $\boldsymbol{u}\in\mathbb{R}^{N\times 3}$. Stage~III uses the predicted field to
    warm-start a classical iterative solver. The prediction lowers the
    iteration count needed to reach the relative residual tolerance
    $\lVert \mathbf{r}\rVert/\lVert \mathbf{F}\rVert \le \mathrm{tol}$.}
  \label{fig:overview}
\end{figure*}

\cref{fig:overview} summarises the FEVessel pipeline. 
The pipeline has three stages. Each vessel
first enters as an unstructured point cloud, with every node
carrying its coordinates together with the local material and load
channels. The pretrained Transolver operator then maps this point
cloud to the full displacement field in a single forward pass. That
prediction finally warm-starts a classical iterative solver,
lowering the iteration count needed to reach the residual tolerance.
The remainder of this section sets out the problem-specific
ingredients of this pipeline. \Cref{subsec:problem} derives the governing
equations and the total potential energy functional that defines the
pretraining loss. \Cref{subsec:encoding} specifies how each vessel is encoded
as the point-cloud input.

\subsection{Problem formulation}\label{subsec:problem}
We consider the static structural analysis of pressure vessels under
internal pressure, governed by the equations of linear
elasticity~\cite{timoshenko1970theory}. Equilibrium requires
\begin{equation}\label{eq:elasticity}
  \nabla \cdot \boldsymbol{\sigma} + \boldsymbol{b} = \boldsymbol{0}
  \quad \text{in } \Omega,
\end{equation}
where $\boldsymbol{\sigma}$ is the Cauchy stress tensor and
$\boldsymbol{b}$ the body force, subject to the boundary conditions
\begin{equation}\label{eq:bc}
  \boldsymbol{u} = \bar{\boldsymbol{u}} \text{ on } \Gamma_u,
  \quad
  \boldsymbol{\sigma} \cdot \boldsymbol{n} = -p\,\boldsymbol{n}
  \text{ on } \Gamma_p,
\end{equation}
where $\Gamma_u$ is the support boundary on which the displacement
$\bar{\boldsymbol{u}}$ is prescribed (fixing the rigid-body modes),
$\Gamma_p$ is the internal pressure surface, $p$ is the internal
pressure, and $\boldsymbol{n}$ is the outward unit normal. The
constitutive law is Hooke's law for isotropic linear elasticity,
\begin{equation}\label{eq:hooke}
  \sigma_{ij}
  = \frac{E}{1 + \nu}\, \varepsilon_{ij}
  + \frac{E\nu}{(1 + \nu)(1 - 2\nu)}\, \varepsilon_{kk}\, \delta_{ij},
\end{equation}
with Young's modulus $E$ and Poisson's ratio $\nu$; here
$\delta_{ij}$ is the Kronecker delta and repeated indices are summed.
The strains follow from the infinitesimal strain--displacement
relation
\begin{equation}\label{eq:strain-disp}
  \varepsilon_{ij}
  = \tfrac{1}{2}\bigl(u_{i,j} + u_{j,i}\bigr),
\end{equation}
giving the strain energy density
\begin{equation}\label{eq:psi}
  \Psi(\boldsymbol{u})
  = \tfrac{1}{2}\, \sigma_{ij}\, \varepsilon_{ij}.
\end{equation}
Neglecting the body force ($\boldsymbol{b} = \boldsymbol{0}$, gravity
being negligible against the internal pressure), the total potential
energy of the vessel is
\begin{equation}\label{eq:Pi-vessel}
  \Pi[\boldsymbol{u}; a]
  = \int_{\Omega(a)} \Psi(\boldsymbol{u})\, \mathrm{d}\Omega
  + \int_{\Gamma_{p}}
      p\, (\boldsymbol{n} \cdot \boldsymbol{u})\, \mathrm{d}\Gamma.
\end{equation}
Substituting $\Pi$ into the energy-form
loss~\cref{eq:pfem-loss-energy} gives the pretraining objective
minimised in this work,
\begin{equation}\label{eq:Lpfem-vessel}
  \mathcal{L}_{\text{PFEM}}^{\text{e}}(\theta)
  = \frac{1}{|\mathcal{D}|}
    \sum_{a \in \mathcal{D}}
    \left[
      \int_{\Omega(a)}
        \Psi\bigl(\mathcal{H}_\theta(a)\bigr)\, \mathrm{d}\Omega
      + \int_{\Gamma_{p}}
        p\, \bigl(\boldsymbol{n} \cdot \mathcal{H}_\theta(a)\bigr)\,
        \mathrm{d}\Gamma
    \right].
\end{equation}
Because the energy form does not enforce the Dirichlet condition on
$\Gamma_u$, we impose it by masking the output component-wise,
\begin{equation}\label{eq:hardBC}
  \boldsymbol{u}(\boldsymbol{x})
  = \boldsymbol{\chi}_{\Gamma_u}(\boldsymbol{x}) \odot
    \mathcal{H}_\theta(a)(\boldsymbol{x}),
\end{equation}
where $\odot$ is the component-wise product and
$\boldsymbol{\chi}_{\Gamma_u}$ is a smooth, distance-based mask whose
entries vanish for the displacement components fixed at the supports,
removing the six rigid-body modes; the support layout is detailed in
\cref{subsec:setup}.

The formulation thus varies along three coupled axes---the
heterogeneous material fields $(E, \nu)$, the geometry $\Omega$, and
the internal pressure $p$---which we call the M/G/B parameter space.
The sampling distributions and ranges are given in
\cref{subsec:setup}.

\subsection{Input encoding for pressure
vessels}\label{subsec:encoding}

Following the input convention of \cref{subsec:transolver}, all
problem information is encoded as pointwise features on an
unstructured point cloud coinciding with the FEM mesh nodes of
$\Omega$, in the form of the Transolver tuple
$\{\boldsymbol{x}^{(i)}, m^{(i)}, b^{(i)}\}_{i=1}^{N}$. For the
pressure vessel, the abstract channels are instantiated as
\begin{equation}\label{eq:encoding}
  m^{(i)} = \bigl(E^{(i)},\, \nu^{(i)}\bigr),
  \qquad
  b^{(i)} = p^{(i)},
\end{equation}
where $E^{(i)}$ and $\nu^{(i)}$ are the local Young's modulus and
Poisson's ratio and $p^{(i)}$ is the prescribed pressure, non-zero on
$\Gamma_p$ and zero elsewhere. A single scalar $p^{(i)}$ thus encodes
both the magnitude and the location of the pressure boundary, so no
separate indicator field is needed. The geometry enters implicitly
through the coordinates $\boldsymbol{x}^{(i)} \in \mathbb{R}^{3}$ (no
explicit geometric descriptor), and the outward normal $\boldsymbol{n}$
is computed from mesh connectivity when the boundary integral of
$\Pi$~\cref{eq:Pi-vessel} is evaluated, rather than learned.

Each point is thus described by its coordinates
$\boldsymbol{x}^{(i)} \in \mathbb{R}^{3}$ and three physical channels
$(E^{(i)}, \nu^{(i)}, p^{(i)})$, and the operator $\mathcal{H}_\theta$
maps this to the three displacement components
$\boldsymbol{u}^{(i)} \in \mathbb{R}^{3}$.

\section{Numerical experiments}\label{sec:experiments}

We evaluate FEVessel on pressure vessel analysis through four sets
of experiments. \cref{subsec:setup} describes the common experimental setup and dataset.
\cref{subsec:mgb} evaluates FEVessel’s generalisation across
material, geometry, and boundary variations. \cref{subsec:warmstart}
quantifies the acceleration obtained by using FEVessel predictions as warm starts for iterative solvers.
\cref{subsec:resolution} demonstrates cross-resolution prediction
from coarse point-cloud input. \cref{subsec:degenerate} demonstrates FEVessel’s robustness 
on two representative types of degenerate meshes, namely inverted and sliver-element meshes,
where FEM fails.

\subsection{Setup and dataset}\label{subsec:setup}

The common experimental setup for all four experiments is summarised
below. All training and inference run on a single NVIDIA GeForce
RTX~5090 GPU, while FEM reference data generation is parallelised
across 8 workers on an AMD Ryzen~9 9950X CPU. All training and
inference times reported in this paper are wall-clock seconds.

We construct a dataset of 3000 pressure vessel configurations by
independently sampling the material, geometry, and boundary parameter axes detailed in \cref{tab:dataset}. The full dataset
is publicly released, with the download link given in the Data
Availability Statement. It is split into $N_{\text{train}}=2500$
training and $N_{\text{test}}=500$ test samples drawn from disjoint
random seeds, ensuring an test set evaluation. 

Throughout the experiments, we evaluate two training regimes of the
same Transolver operator, the physics-informed Transolver and the
supervised Transolver, which share the same network architecture and
differ only in the training supervision signal. The physics-informed
Transolver is the operator that FEVessel pretrains and the standard
configuration adopted in this work, trained on the total potential
energy functional alone without any displacement labels. The
supervised Transolver is an ablation variant that replaces the energy
loss with a supervised relative $L^2$ loss against the FEM reference
displacements. All other architecture and training settings remain
identical to the physics-informed Transolver. Because the
physics-informed Transolver trains on the minimum potential energy
alone without displacement labels, the train/test partition raises no
label-leakage concern for the model. The partition is retained so
that the supervised baselines, the supervised Transolver and FNO, are
evaluated on the test set. The physics-informed Transolver itself does
not require a separate data split.

\begin{figure*}[htbp]
  \centering
  \includegraphics[width=\textwidth]{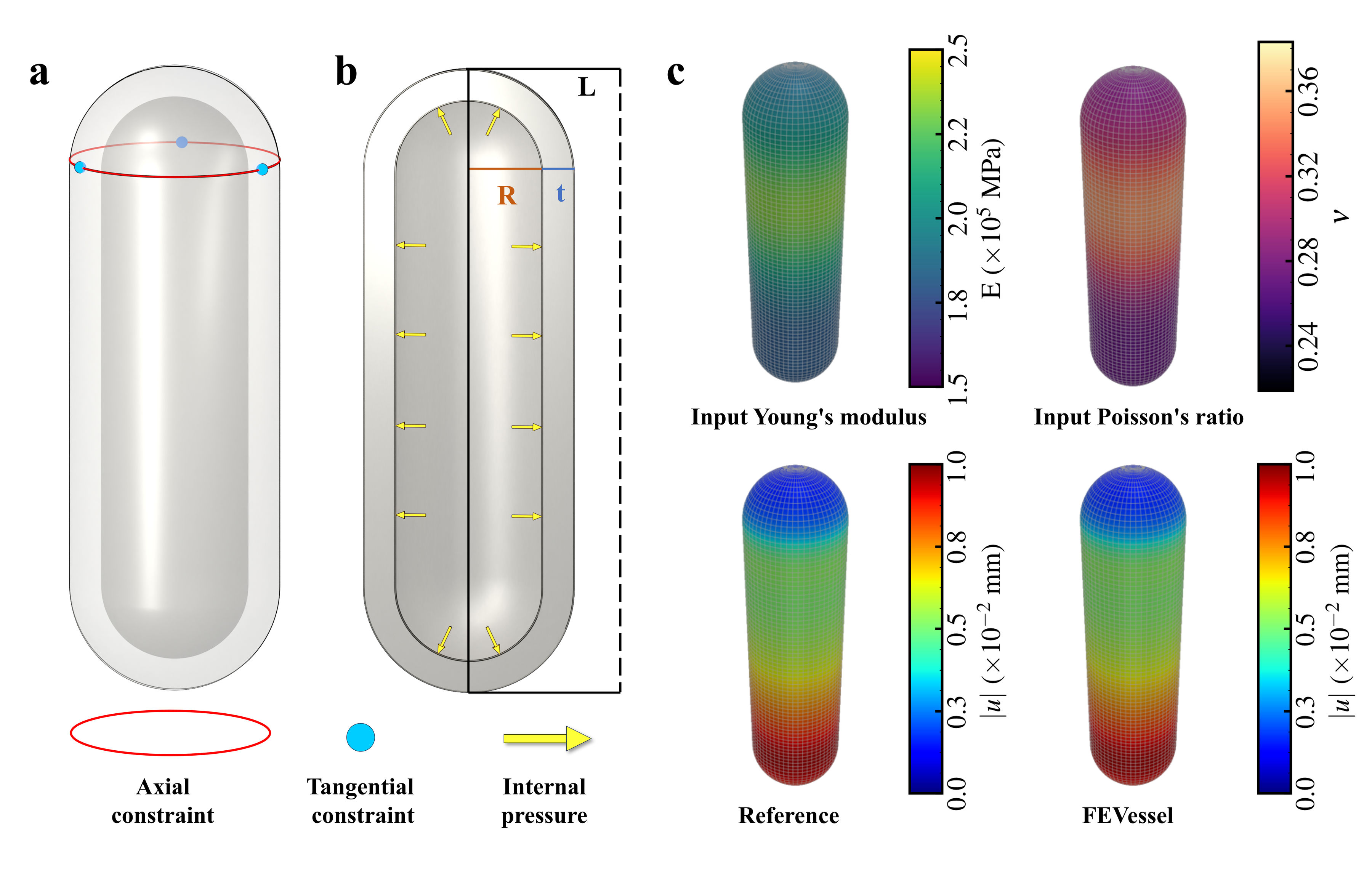}
  \caption{Problem set-up and a representative FEVessel solution.
  (a)~Pressure vessel geometry with the axial and tangential supports that
  remove all six rigid-body modes. (b)~Geometric parameters and the internal pressure load,
  with inner radius $R$, wall thickness $t$, and cylinder
  length $L$. (c)~One heterogeneous sample, with the input $E$ and $\nu$ fields (top row) and the FEM reference
  and FEVessel displacement magnitudes $|\boldsymbol{u}|$ compared on a shared
  colour scale (bottom row).}
  \label{fig:setup}
\end{figure*}

\Cref{fig:setup} gives the experimental setup.
Because a uniform internal pressure on a closed surface is
self-equilibrated, the six rigid-body modes must be fixed
explicitly.
An axial constraint along the whole circumferential ring
(red in \cref{fig:setup}a) and tangential constraints at
three equally spaced points on that ring (blue dots) together
remove all six rigid-body modes while leaving the
pressure-induced deformation unrestrained.
We then specify the geometry. To represent a realistic pressure
vessel, \cref{fig:setup}b adopts a cylinder closed by two
hemispherical heads, parameterised by inner radius $R$,
wall-thickness ratio $t/R$, and cylinder length ratio $L/R$, each
sampled uniformly from the ranges in \cref{tab:dataset}. The
internal pressure $p$ is drawn uniformly from $[50, 200]$~MPa and
applied on the inner surface, shown as the yellow arrows. Finally,
we set the material. The fields are spatially heterogeneous, with
the Young's modulus $E$ and Poisson's ratio $\nu$ each realised as
an independent Gaussian random field (GRF) on a $24\times24\times64$
background grid and interpolated to the FEM mesh nodes, with mean
values $\bar{E}=200$~GPa and $\bar{\nu}=0.30$, standard deviations
$\sigma_E=20$~GPa and $\sigma_\nu=0.03$, clipped to the physically
admissible ranges $E\in[150,250]$~GPa and $\nu\in[0.20,0.40]$.
\Cref{fig:setup}c shows one representative sample, with its sampled
$E$ and $\nu$ fields in the top row and the FEM reference and FEVessel
displacement magnitudes compared on a shared colour scale in the
bottom row, where the two fields show close visual agreement.

\begin{table}[!t]
  \centering
  \caption{Parameter ranges and sampling distributions for the
    pressure vessel dataset. In the Distribution column, GRF
    denotes a clipped Gaussian random field sampled per node and
    interpolated to the mesh, whereas Uniform denotes independent
    uniform sampling of a scalar parameter.}
  \label{tab:dataset}
  \begin{tabular*}{\linewidth}{@{\extracolsep{\fill}}lll@{}}
    \toprule
    Parameter & Range & Distribution \\
    \midrule
    Young's modulus $\bar{E}$    & $200\pm20$~GPa, clipped $[150,250]$    & GRF \\
    Poisson's ratio $\bar{\nu}$  & $0.30\pm0.03$, clipped $[0.20,0.40]$   & GRF \\
    Inner radius $R$             & $[2.5, 6.5]$~mm                        & Uniform \\
    Thickness ratio $t/R$        & $[0.12, 0.35]$                         & Uniform \\
    Length ratio $L/R$           & $[7.0, 14.0]$                          & Uniform \\
    Internal pressure $p$        & $[50, 200]$~MPa                        & Uniform \\
    \bottomrule
  \end{tabular*}
\end{table}

High-fidelity reference solutions are generated by
linear elasticity solver using C3D8R elements (8-node hexahedral,
reduced integration) with hourglass stabilisation factor
$\alpha=0.2$. Mesh density is adaptively scaled to a target element
size of approximately $0.5$~mm in all directions, with
circumferential divisions rounded to multiples of 4 to preserve
cylindrical symmetry. Force-balance quality is verified per sample
by requiring the reaction-to-load ratio to remain below $10^{-5}$,
and all cases pass mesh convergence studies. These reference
solutions are used only to evaluate the trained operator on the test set. They are never used as supervision during FEVessel
training, which relies exclusively on the physics-informed energy
loss.

For all experiments we report the relative $L^2$ error of the
displacement vector field on the test set,
\begin{equation}\label{eq:l2error}
  \varepsilon_{L^2}(\boldsymbol{u}) =
  \sqrt{
    \dfrac{\displaystyle\int_\Omega
    \left\|
      \boldsymbol{u}_{\text{pred}}(\boldsymbol{x})
      - \boldsymbol{u}_{\text{ref}}(\boldsymbol{x})
    \right\|_2^2 \,\mathrm{d}\Omega}
    {\displaystyle\int_\Omega
    \left\|
      \boldsymbol{u}_{\text{ref}}(\boldsymbol{x})
    \right\|_2^2 \,\mathrm{d}\Omega}
  },
\end{equation}
where $\boldsymbol{u}_{\text{pred}}$ is the operator prediction,
$\boldsymbol{u}_{\text{ref}}$ the FEM reference, and $\|\cdot\|_2$ the
Euclidean norm over the three displacement components at each point.
We additionally report the relative $H^1$ semi-norm error, which
measures the accuracy of the displacement gradient---equivalently the
strain, and through Hooke's law the stress---rather than of the
displacement itself,
\begin{equation}\label{eq:h1error}
  \varepsilon_{H^1}(\boldsymbol{u}) =
  \sqrt{
    \dfrac{\displaystyle\int_\Omega
    \left\| \nabla\boldsymbol{u}_{\text{pred}}(\boldsymbol{x})
      - \nabla\boldsymbol{u}_{\text{ref}}(\boldsymbol{x}) \right\|_F^2
    \,\mathrm{d}\Omega}
    {\displaystyle\int_\Omega
    \left\| \nabla\boldsymbol{u}_{\text{ref}}(\boldsymbol{x}) \right\|_F^2
    \,\mathrm{d}\Omega}
  },
\end{equation}
where $\nabla\boldsymbol{u}$ is the displacement gradient and
$\|\cdot\|_F$ the Frobenius norm, evaluated per element with the same
C3D8R shape-function derivatives used by the energy loss. A prediction
can match the displacement field well (small $\varepsilon_{L^2}$) yet
have an inaccurate gradient (large $\varepsilon_{H^1}$), so
$\varepsilon_{H^1}$ is the more stringent, engineering-relevant metric
for the stress field.
The error is computed on the full displacement vector $\boldsymbol{u}$
(not the scalar magnitude $|\boldsymbol{u}|$), and statistics are
reported as mean $\pm$ standard deviation over the test set.
Unless stated otherwise, a single reported error is the converged
value, taken as the lowest test error over the final third of
training, by which point every model has reached its plateau.

\subsection{Generalisation across material, geometry, and
boundary variations}\label{subsec:mgb}

This experiment evaluates the generalisation capability of FEVessel.
We compare the physics-informed Transolver that FEVessel pretrains with
two baselines on the dataset of \cref{subsec:setup}: a supervised
Transolver variant that replaces the energy loss with label
supervision, and a structured-grid FNO.
The input follows the encoding defined in \cref{subsec:encoding}. The pressure
vessel index $i$ ranges over the $3000$ dataset configurations of
\cref{subsec:setup}, and the point index $j$ ranges over the mesh
nodes of vessel $i$, so that $j = 1, \dots, N^{(i)}$. The node count
$N^{(i)}$ is vessel-specific, because each geometry is meshed
independently and therefore yields a different number of nodes. For
vessel $i$ and point $j$, each input point carries its spatial
coordinate $\boldsymbol{x}_j^{(i)}$, the local material pair
$(E_j^{(i)}, \nu_j^{(i)})$, and the pressure channel $b_j^{(i)}$,
and the operator returns the displacement field
$\{ \boldsymbol{u}_j^{(i)} \}_{j=1}^{N^{(i)}}$ on the same points.

We train the physics-informed Transolver on the full M/G/B training
distribution and evaluate it on the 500 test cases in a single
forward pass with no fine-tuning. The operator is a
4-layer Transolver backbone with 128 hidden dimensions, 8 attention
heads, and 32 physics-aware tokens, trained with the Adam optimiser
at initial learning rate $\eta = 2\times10^{-3}$ and cosine
annealing to $\eta_{\min} = 1\times10^{-5}$ over 1300 epochs, with
batch size~4 and gradient clipping at norm~1.0. To separate the
contribution of the point-cloud representation from the effect of
the physics-informed objective, we train a supervised Transolver variant
with the same architecture and training settings but replace the
energy loss with the relative $L^2$ error against the FEM
displacement labels. The physics-informed and supervised Transolver models take 29.4~h and 25.6~h to train, respectively.

For comparison, we train a supervised FNO baseline on the same
dataset. Because FNO requires a structured input, each pressure
vessel is resampled onto a body-fitted shell grid of resolution
$N_s \times N_\theta \times N_\xi = 128\times64\times8$, where $s$
is the meridian arc-length coordinate running from pole to pole,
$\theta$ the circumferential angle, and $\xi$ the through-thickness
coordinate, so that every grid point lies inside the vessel wall.
The material fields $(E, \nu)$ are interpolated to the grid nodes,
and the input further carries the pressure load $p$ and the
geometry scalars $(R, t/R, L/R)$. The FNO has 4 layers with 64
hidden channels and 16 retained Fourier modes per spectral
direction, applied on the $(s, \theta)$ mid-surface with the $\xi$
layers folded into the channel dimension. It returns the
three-component displacement field $(u_x, u_y, u_z)$ on the grid,
mapped back to the mesh nodes by trilinear interpolation in
$(s, \theta, \xi)$ parameter space, and is trained with a
mean-squared error (MSE) against the FEM reference displacements.
All three models are trained for 1300 epochs with the same Adam
optimiser and cosine learning-rate schedule.

\begin{table*}[t]
  \centering
  \caption{Comparison of FNO, the physics-informed Transolver, and the
  supervised Transolver on the M/G/B pressure vessel test set (500 cases).
  FNO uses a body-fitted structured shell grid, whereas both Transolver
  variants operate on the unstructured point cloud. All
  models are trained for 1300 epochs with the Adam optimiser and
  cosine learning-rate annealing, by which point every test curve
  has reached its plateau. The reported $L^2$ (displacement) and
  $H^1$ semi-norm (strain) errors are per-sample means over the 500
  test cases at the converged model (\cref{eq:l2error,eq:h1error}).
  Data generation is the
  wall-clock time to produce the 2500 FEM training labels on 8 CPU
  workers, and the physics-informed Transolver is label-free and needs none.
  The inference column reports the mean single-sample forward-inference
  wall-clock time over the 500 test cases (\emph{model}), against the mean
  single-sample classical FEM solve time (\emph{FEM}, AMG-preconditioned CG
  at the $10^{-3}$ engineering tolerance from a zero initial guess); both are
  means over the 500 test cases at batch size one, and the FNO \emph{model}
  figure includes the mandatory resampling of the mesh onto the structured
  shell grid and the mapping of the prediction back to the mesh nodes.}
  \label{tab:fno_pfem_comparison}
  \setlength{\tabcolsep}{2pt}%
  \begin{tabular*}{\textwidth}{@{\extracolsep{\fill}}lllccccc@{}}
    \toprule
    Method & Loss (Supervision) & Hidden / Modes & Rel $L^2(\boldsymbol{u})$ & Rel $H^1$ semi-norm & Data generation & Train time & \begin{tabular}[c]{@{}c@{}}Inference (s)\\ model / FEM\end{tabular} \\
    \midrule
    FNO & MSE (supervised) & 64 / 16 modes & 0.57\% & 9.72\% & 6.6 h & 18.3 h & 0.042 / 2.30 \\
    Physics-informed Transolver & Energy $\Pi$ (label-free) & 128 / 32 tokens & 1.35\% & 2.07\% & 0 h & 29.4 h & 0.0130 / 2.30 \\
    Supervised Transolver & Relative $L^2$ (supervised) & 128 / 32 tokens & 0.14\% & 1.70\% & 6.6 h & 25.6 h & 0.0130 / 2.30 \\
    \bottomrule
  \end{tabular*}
\end{table*}

\begin{figure*}[!t]
  \centering
  \includegraphics[width=\linewidth]{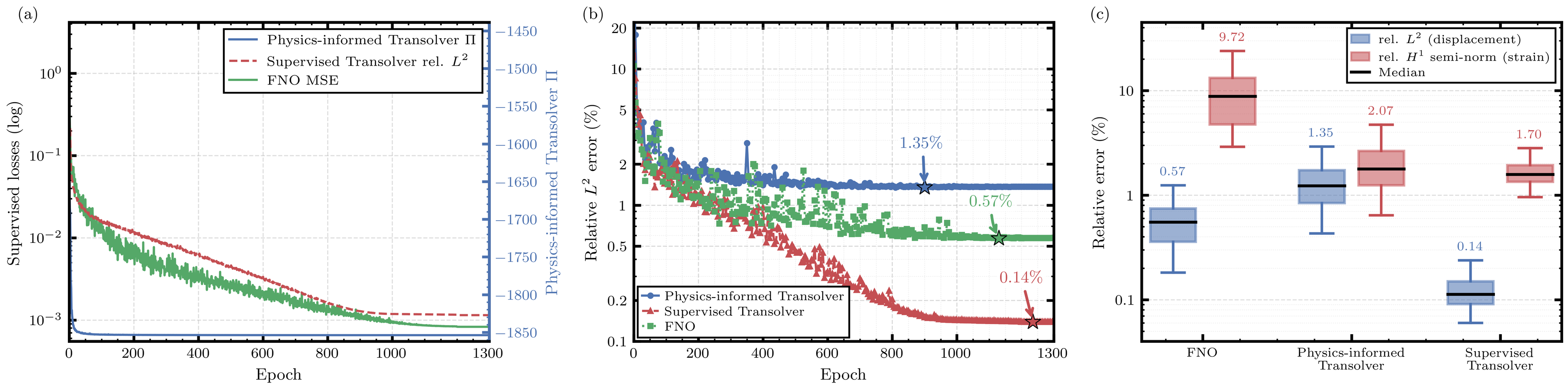}
  \caption{Training dynamics of FNO, the physics-informed Transolver,
    and the supervised Transolver. (a)~The left logarithmic axis carries the two
    supervised losses, the FNO MSE and the supervised Transolver
    relative $L^2$ loss. The right linear axis carries the
    physics-informed Transolver energy $\Pi$, which is minimised towards a negative
    value and so cannot share the scale of the supervised losses.
    (b)~Test relative $L^2$ error of $\boldsymbol{u}$ versus epoch.
    The supervised Transolver converges to $0.14\%$, FNO to $0.57\%$,
    and the physics-informed Transolver to $1.35\%$.
    (c)~Per-sample distributions of the relative $L^2$ (displacement)
    and $H^1$ semi-norm (strain) errors over the 500 test cases, on a
    logarithmic axis (boxes: interquartile range and median).}
  \label{fig:training-curves}
\end{figure*}

\Cref{fig:training-curves}a shows that all three training
objectives, the physics-informed Transolver energy $\Pi$, the FNO MSE, and the
supervised Transolver relative $L^2$ loss, converge and flatten within
the 1300 epochs, confirming that the models are adequately trained.
On the displacement itself all three models are accurate
(\cref{tab:fno_pfem_comparison} and \cref{fig:training-curves}b):
the supervised Transolver converges to a $0.14\%$ relative
displacement-vector $L^2$ error, FNO to $0.57\%$, and the
physics-informed Transolver to $1.35\%$. The strain tells a
different story. The $H^1$ semi-norm error reaches $9.72\%$ for FNO
but only $2.07\%$ and $1.70\%$ for the physics-informed and
supervised Transolver (\cref{fig:training-curves}c), and this is
the metric that carries the stress and therefore the engineering
assessment.

The strain gap stems from the structured representation rather than
the network or the supervision signal. FNO computes on a fixed
grid, so the fields must cross the grid--mesh interface twice: the
inputs are resampled from the mesh onto the grid, and the
prediction is interpolated back to the mesh nodes for evaluation.
Projecting the FEM reference displacement itself through these two
interpolations in succession, from the mesh onto the grid and then
back to the mesh nodes, with no learning involved, already incurs a
$0.64\%$ $L^2$ and an $11.8\%$ $H^1$ error on the test set. Both FNO errors
therefore coincide with the transfer error of the representation
itself, not with a limit of optimisation: the displacement, a
smooth field, survives the transfer almost intact, but its gradient
concentrates through the thin wall, no fixed grid--mesh
interpolation preserves it, and differentiation then amplifies the
loss (\cref{subsec:geometric-representations}). The point-cloud
Transolver computes directly on the mesh nodes, so no such transfer
exists and the strain stays accurate at $2.07\%$, about
$4.7\times$ better than FNO.

The inference-time gap is pronounced. Once trained, the
physics-informed Transolver predicts a new vessel in $0.0130$~s,
against $2.30$~s for the FEM solve at the same $10^{-3}$ tolerance,
about $177\times$ faster; FNO needs $0.042$~s, dominated by the
mandatory mesh--grid resampling (\cref{tab:fno_pfem_comparison}).
The cost is training time. The physics-informed Transolver trains
longer than FNO, $29.4$~h against $18.3$~h, and on the displacement
$L^2$ error the two supervised models are ahead ($0.14\%$ and
$0.57\%$ against $1.35\%$). But both supervised models depend on
the 2500 precomputed FEM solutions, which take $6.6$~h to generate
here and are often expensive or simply unavailable in practice, and
on the engineering-relevant strain metric the ordering reverses,
with the physics-informed Transolver about $4.7\times$ more
accurate than FNO ($2.07\%$ against $9.72\%$). The physics-informed
Transolver is therefore the more deployable choice, delivering
engineering-accurate stress with no labels and no mesh--grid
transfer.

\subsection{Warm-start acceleration}\label{subsec:warmstart}

This experiment quantifies the acceleration obtained by using the
FEVessel prediction as the initial guess for a classical iterative
solver, while preserving engineering-grade accuracy. For each of the
500 test cases we assemble the stiffness matrix $\mathbf{K}$ and load
vector $\mathbf{F}$ and solve $\mathbf{K}\mathbf{U}=\mathbf{F}$ with three solvers, 
unpreconditioned conjugate gradient (CG)~\cite{saad2003iterative},
algebraic multigrid (AMG)~\cite{vanek1996algebraic}, and an
AMG-preconditioned CG whose smoothed-aggregation hierarchy is
supplied with the six elastic rigid-body near-null-space modes
(AMG-PCG). The zero initial guess $\mathbf{U}^{(0)}=\mathbf{0}$ is
compared against the FEVessel warm-start
$\mathbf{U}^{(0)}=\mathcal{H}_\theta(a)$ across three relative
residual tolerances $\text{tol}\in\{10^{-2},10^{-3},10^{-4}\}$, with
convergence declared when $\|\mathbf{r}\|/\|\mathbf{F}\|$ falls below
$\text{tol}$.

\begin{figure*}[!t]
  \centering
  \includegraphics[width=\linewidth]{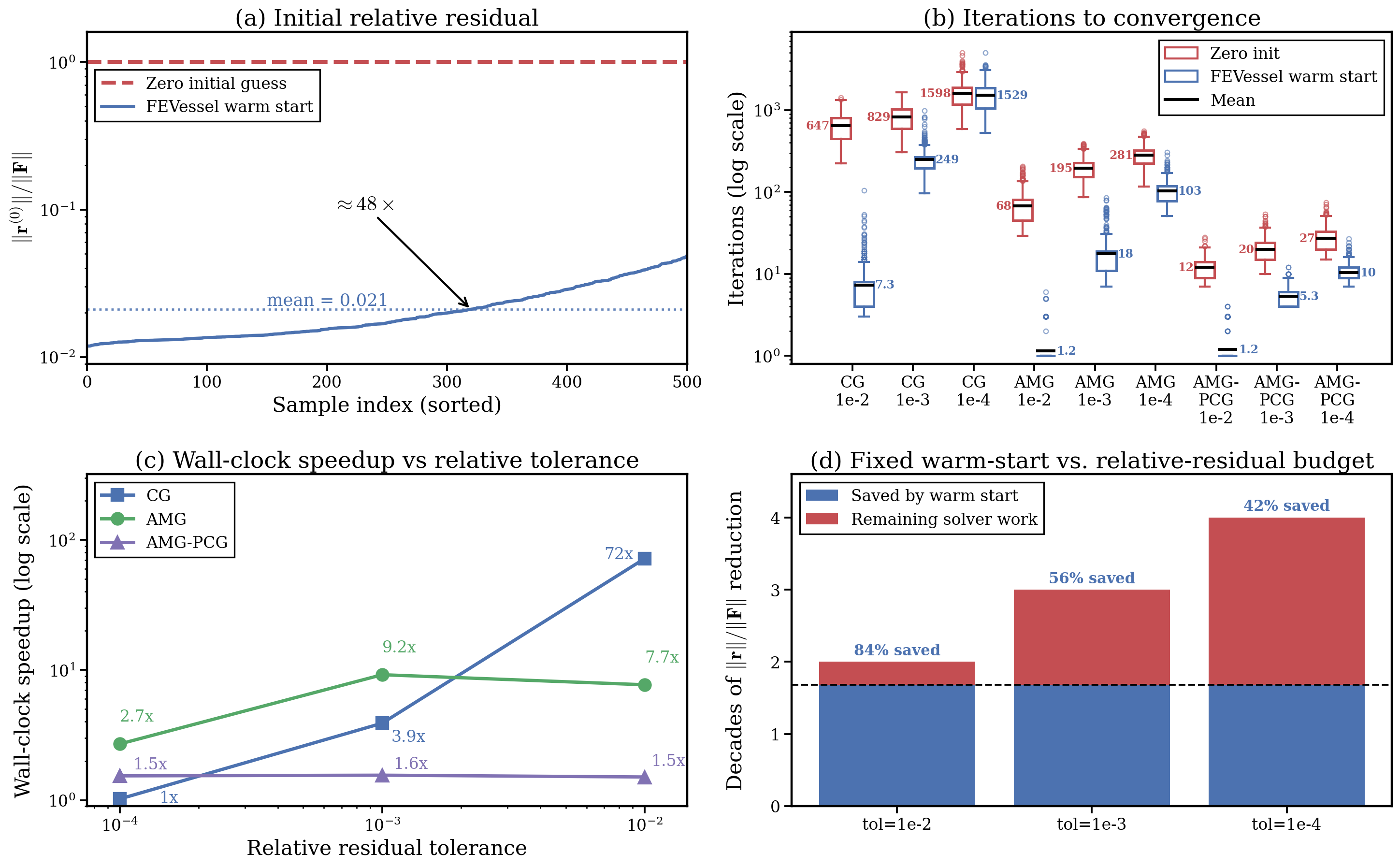}
  \caption{Warm-start behaviour over 500 test cases.
    (a)~Sorted initial relative residual
    $\|\mathbf{r}_0\|/\|\mathbf{F}\|$, reduced from $1$ to about $0.021$, a
    factor of roughly $48$, by the warm-start. (b)~Iterations
    to convergence (CG, AMG, and AMG-PCG). (c)~Wall-clock speedup
    versus relative residual tolerance. (d)~Decomposition of the
    relative-residual speedup ceiling.}
  \label{fig:warmstart-conv}
\end{figure*}

Before any solver iteration, FEVessel already places the solution close
to convergence. \Cref{fig:warmstart-conv}a shows the initial
relative residual sorted over the 500 test cases. A zero initial
guess starts at $\|\mathbf{r}_0\|/\|\mathbf{F}\|=1$ by construction, whereas
the FEVessel warm-start starts at $0.021$ on average, a reduction of
about $48\times$.

The reduced initial residual translates directly into fewer
iterations. \Cref{tab:warmstart} reports the mean iteration count
over the 500 test cases at each tolerance, and all three solvers
drop markedly. At the commonly used tolerance of $10^{-3}$, AMG
falls most steeply, from 195 to 17.6, with CG from 829 to 249 and
AMG-PCG from 19.89 to 5.33, giving mean iteration speedups of
$11.0\times$, $3.3\times$, and $3.73\times$. The boxplots in
\cref{fig:warmstart-conv}b confirm that this reduction holds across
the test set rather than being driven by a few favourable samples.

\begin{table}[!t]
  \caption{Iterative-solver convergence with zero versus FEVessel
    warm-start initialisation, averaged over 500 test cases. The
    iteration speedup is the ratio of the mean iteration counts
    $\bar{n}_{\text{zero}}/\bar{n}_{\text{FEVessel}}$, and the
    wall-clock speedup is likewise the ratio of the mean end-to-end
    times, including the FEVessel forward pass and, for AMG-based
    methods, the multigrid hierarchy and preconditioner setup. Both
    ratios are computed from the unrounded means; because the
    per-iteration cost varies with the vessel-specific mesh size,
    the two ratios need not coincide.}
  \label{tab:warmstart}
  \begin{tabular*}{\linewidth}{@{\extracolsep{\fill}}llcccc@{}}
    \toprule
    \multirow{2}{*}{Solver} & \multirow{2}{*}{\shortstack{Rel.\\tol.}} &
    \multicolumn{2}{c}{Iterations (mean)} &
    \multirow{2}{*}{\shortstack{Iter.\\speedup}} &
    \multirow{2}{*}{\shortstack{Wall-clock\\speedup}} \\
    \cmidrule(lr){3-4}
                            &           & Zero & FEVessel &              &  \\
    \midrule
    CG  & $10^{-2}$ & 647  & 7.35 & $88.0\times$ & $71.9\times$ \\
    CG  & $10^{-3}$ & 829  & 249  & $3.3\times$  & $3.9\times$  \\
    CG  & $10^{-4}$ & 1598 & 1529 & $1.05\times$ & $1.02\times$          \\
    \midrule
    AMG & $10^{-2}$ & 67.5 & 1.16 & $58.4\times$ & $7.7\times$  \\
    AMG & $10^{-3}$ & 195  & 17.6 & $11.0\times$ & $9.2\times$  \\
    AMG & $10^{-4}$ & 281  & 103  & $2.7\times$  & $2.7\times$           \\
    \midrule
    AMG-PCG & $10^{-2}$ & 12.02 & 1.20  & $10.02\times$ & $1.51\times$ \\
    AMG-PCG & $10^{-3}$ & 19.89 & 5.33  & $3.73\times$  & $1.55\times$ \\
    AMG-PCG & $10^{-4}$ & 27.36 & 10.43 & $2.62\times$  & $1.53\times$ \\
    \bottomrule
  \end{tabular*}
\end{table}

\begin{figure*}[!t]
  \centering
  \includegraphics[width=\linewidth]%
    {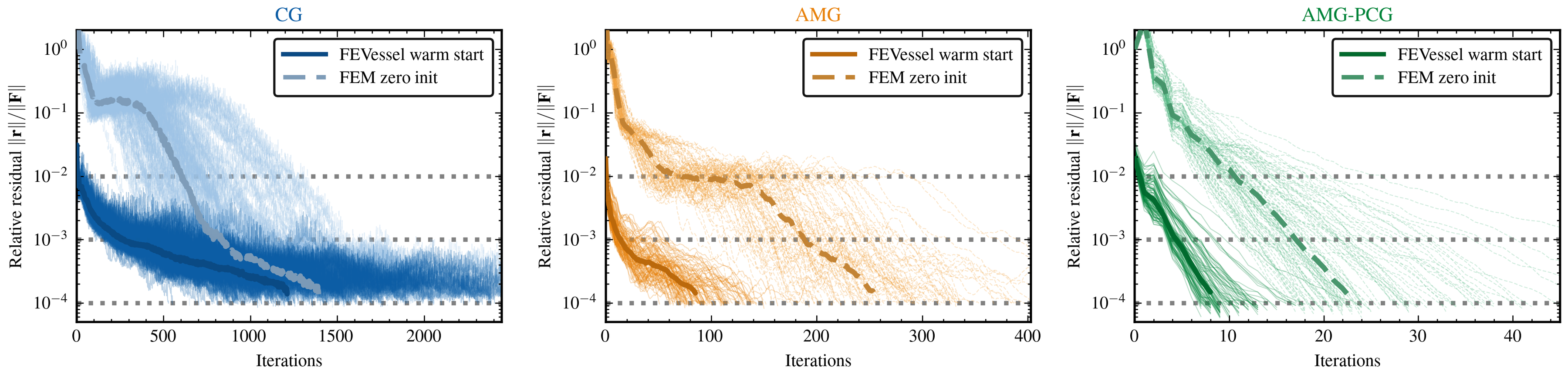}
  \caption{Convergence histories of the relative residual
    $\|\mathbf{r}\|/\|\mathbf{F}\|$ for the warm-start study, shown
    separately for CG, AMG, and AMG-PCG. Within each panel the solid lines denote the FEVessel warm
    start and the dashed lines the FEM zero initial guess,
    with the thick line the median over 100 test cases, drawn at
    random from the 500-case test set and shown as thin lines. The dotted
    horizontal lines mark the convergence tolerances $10^{-2}$,
    $10^{-3}$, and $10^{-4}$.}
  \label{fig:warmstart-history}
\end{figure*}

The reduction in iteration count comes from the lower residual
starting point, not from a faster convergence rate. The relative
residual convergence histories of the three solvers are shown in
\cref{fig:warmstart-history}. The convergence rate is set by the
solver and its preconditioner and is independent of the initial
guess. For CG, AMG, and AMG-PCG alike, the warm-start cluster
(solid lines) lies below the zero-initial-guess cluster (dashed
lines), starting lower, crossing each tolerance line earlier, and
descending at the same slope.

How these iteration savings translate into wall-clock time differs
across solvers (\cref{fig:warmstart-conv}c), because the fixed time
costs around each solver do not change with the initial guess. CG
carries no setup cost, so its wall-clock speedup closely follows the
iteration savings, falling monotonically from $71.9\times$
at $10^{-2}$ to $3.9\times$ at $10^{-3}$ and $1.02\times$ at
$10^{-4}$. AMG behaves differently, because its multigrid hierarchy
must be rebuilt for every distinct stiffness matrix, a setup cost
that is independent of the initial guess and that the warm start
cannot remove. \Cref{app:warmstart-walltime} shows the end-to-end
speedup is governed jointly by this setup cost and the iteration
saving (\cref{eq:app-speedup}). At $10^{-2}$ the setup cost holds
the speedup to $7.7\times$ despite a $58.4\times$ iteration speedup,
after which it peaks at $9.2\times$ at the engineering tolerance
$10^{-3}$ and recedes to $2.7\times$ at $10^{-4}$.

This setup-dominated effect is more pronounced for AMG-PCG.
AMG-PCG reaches the three tolerances in only 12.02, 19.89, and
27.36 zero-start iterations on average, which the FEVessel warm start
further lowers to 1.20, 5.33, and 10.43. Because the hierarchy
setup is paid under both initialisations and the near-null-space
preconditioner already keeps the zero-start count low, the
iteration saving $\Delta n$ is small. The end-to-end gain therefore
stays smaller and nearly flat, between $1.51\times$ and $1.55\times$
across the three tolerances.

The weakening of the speedup as the tolerance tightens is
structural, and \cref{app:warmstart-theory} traces it to the fixed
residual saving of the warm start. The warm start lowers the
initial relative residual by about 1.68 decades
($\eta_0\approx0.021$). This saving carries no tolerance, so its
share of the total residual budget falls as the tolerance tightens.
The decomposition in \cref{fig:warmstart-conv}d shows this share
dropping from about $84\%$ at $10^{-2}$ to $56\%$ at $10^{-3}$ and
$42\%$ at $10^{-4}$, so the speedup weakens accordingly. Whatever
the tolerance, once the solver converges the warm start returns the
same solution of the same linear system as the zero guess at the
prescribed tolerance, so engineering-grade accuracy is preserved
exactly.

\subsection{Cross-resolution prediction}\label{subsec:resolution}

This experiment evaluates FEVessel's ability to keep its accuracy
as the input point density varies across pressure vessel
meshes. We reuse the FEVessel model trained in \cref{subsec:mgb}
without any fine-tuning. At inference time we subsample each
test vessel's input point cloud at densities
$\rho \in \{1\%, 2\%, \dots, 100\%\}$, pass the subsampled
cloud through the operator, and reconstruct the displacement
field at the vessel's full set of mesh-node query locations by
$k$-nearest-neighbour interpolation ($k = 4$). The
reconstructed field is compared against the high-fidelity FEM
reference. Statistical curves are averaged over $50$ test
vessels and $3$ random subsampling seeds. For the single-vessel
examination we select one specific representative vessel,
Sample~ID~0334, whose mesh contains
$N_{\text{fine}} = 41{,}768$ nodes, so its subsampled inputs
span $\sim\!418$ to $41{,}768$ points.

\begin{figure*}[!t]
  \centering
  \includegraphics[width=\linewidth]{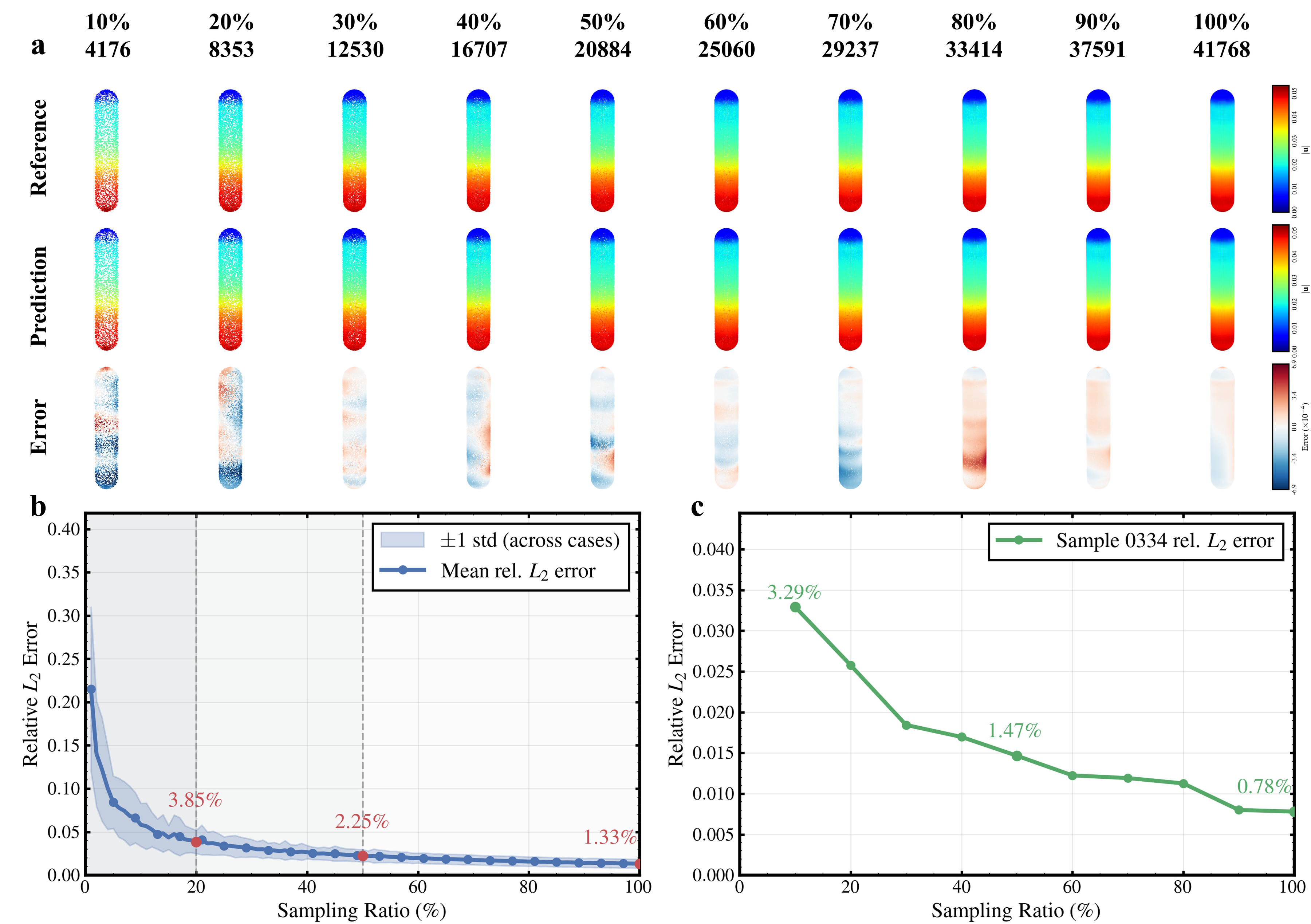}
  \caption{Cross-resolution prediction by FEVessel.
    (a)~Displacement magnitude $|\boldsymbol{u}|$ on Sample~ID~0334 at
    ten input densities from $10\%$ to $100\%$ (FEM reference, FEVessel
    prediction, and vector error). (b)~Mean relative $L^2$ error of
    $\boldsymbol{u}$ versus input ratio over $50$ test vessels, with the
    undersampling, intermediate, and saturation regimes shaded and
    the operator-accuracy floor marked. (c)~Per-sample relative
    $L^2$ error of $\boldsymbol{u}$ versus input ratio for
    Sample~ID~0334, the case shown in~(a).}
  \label{fig:resolution}
\end{figure*}

On this case, FEVessel predicts accurately at every density
tested, and the accuracy improves monotonically as the input
density rises. Training presents each vessel only at its full
mesh, and these coarsened densities never occur during
training, so the retained accuracy reflects genuine
cross-resolution generalisation rather than memorisation of a
particular grid. \Cref{fig:resolution}a shows the prediction
at ten densities from $10\%$ to $100\%$, each reproducing the
smooth displacement field of the FEM reference, and even at
$10\%$ input ($4{,}176$ points) the field stays qualitatively
correct. The per-sample relative $L^2$ error in
\cref{fig:resolution}c falls from $3.29\%$ at $10\%$ to
$1.47\%$ at $50\%$ and $0.78\%$ at $100\%$.

This trend is not limited to the single case. The statistical
sweep across $50$ test vessels in \cref{fig:resolution}b shows
three regimes. In the undersampling regime
($\rho \lesssim 20\%$) the mean vector-field error falls
rapidly, from $21.5\% \pm 9.5\%$ at $\rho = 1\%$ to
$3.85\% \pm 1.4\%$ at $\rho = 20\%$, and at the lowest
densities a non-trivial share of this error comes from the
$k$-nearest-neighbour reconstruction rather than the operator,
since at $\rho = 1\%$ only about one point in a hundred ($418$
of $41{,}768$ for the representative vessel) is interpolated to
the full set of query locations. In the intermediate regime
($20\% \lesssim \rho \lesssim 50\%$) the decline slows,
reaching $2.25\% \pm 0.7\%$ at $\rho = 50\%$. Beyond
$\rho \approx 50\%$ the error saturates, falling to
$1.33\% \pm 0.5\%$ at full resolution. The standard deviation
across vessels falls from $9.5\%$ to $0.5\%$ over the same
range, so FEVessel's accuracy becomes not only higher but more
uniform across the test set as the input density increases.

This curve, which falls fast and then levels off, gives a
direct engineering choice. The operating point at
$\rho \approx 30\%$ ($\sim\!12{,}500$ points for the
representative vessel, mean error $3.12\%$) is a favourable
trade-off, reducing the input to
about a third of the fully resolved mesh while keeping the
error within $2.3\times$ of the full-resolution baseline
($1.33\%$). For early-stage design exploration, where rapid
evaluation matters more than verification-grade accuracy, FEVessel
can therefore run on coarse industrial meshes without
resolution-specific retraining and without substantial loss of
fidelity.

\subsection{Mesh independence on degenerate
meshes}\label{subsec:degenerate}

Automated CAD-based mesh generation routinely produces degenerate
meshes that defeat conventional FEM. Whether FEVessel still solves on
these meshes is the most distinctive capability we evaluate. To
test this, we start from a high-quality reference mesh on which FEM
converges ($31{,}536$ hexahedral elements) and degenerate a
controlled fraction of the elements in two ways. Inverted elements
reverse the bottom-face winding so that the isoparametric Jacobian
determinant turns negative while the node coordinates stay fixed.
Sliver elements translate one node by $90\%$ of its distance to the
element centroid, severely flattening the element without inverting
it. The degeneration is applied to
$\rho_d \in \{1\%, 5\%, 10\%, 25\%, 50\%, 75\%\}$ of the elements,
chosen uniformly at random, and the results are averaged over $30$
test vessels. Both solvers run on the degenerate mesh and are
compared against the FEM reference on the original mesh.

\begin{figure*}[!t]
  \centering
  \includegraphics[width=\linewidth]{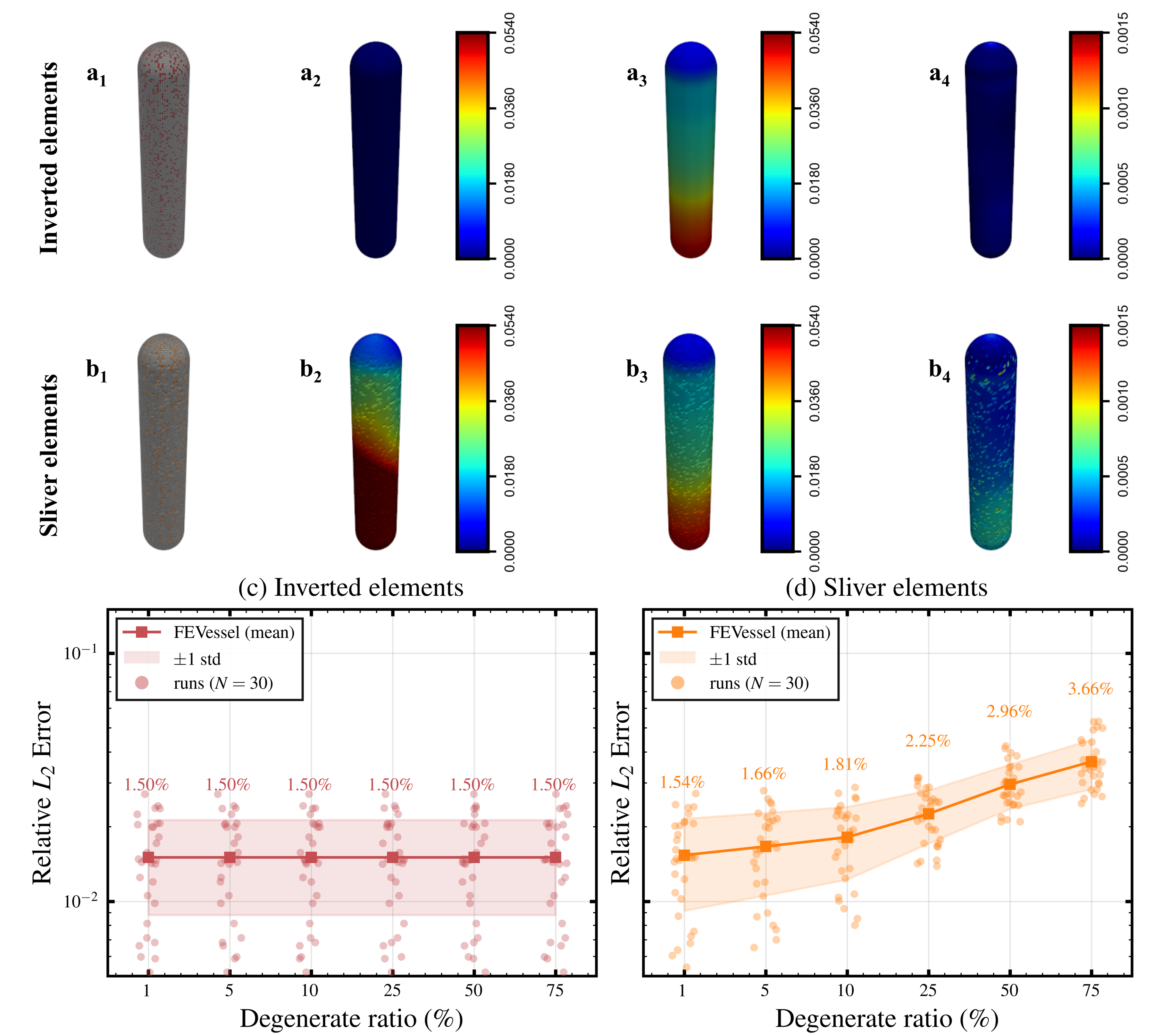}
  \caption{FEM versus FEVessel on degenerate meshes at $\rho_d = 25\%$.
    Rows~a (inverted) and~b (sliver) share four columns:
    (a1,\,b1)~the mesh with the degenerate elements highlighted,
    (a2,\,b2)~the forced-FEM diagnostic $|\boldsymbol{u}|$,
    (a3,\,b3)~the FEVessel $|\boldsymbol{u}|$, and (a4,\,b4)~the
    FEVessel error on a finer colour scale. (c,\,d)~Mean relative
    $L^2$ error of $\boldsymbol{u}$ versus degeneration ratio over
    $30$ vessels for the inverted and sliver cases, from $1.50\%$~(c)
    to $3.66\%$~(d).}
  \label{fig:degenerate}
\end{figure*}

On both degenerate mesh types, conventional FEM fails on all twelve
configurations (six ratios $\times$ two types), but the two failure
mechanisms differ sharply. \Cref{fig:degenerate} contrasts them at
$\rho_d = 25\%$. Rows~a (inverted) and~b (sliver) show, from left
to right, the degenerate mesh, the forced-FEM diagnostic
$|\boldsymbol{u}|$, the FEVessel prediction $|\boldsymbol{u}|$, and the
FEVessel error, and panels~(c) and~(d) plot the mean error against the
degeneration ratio for the two cases. The two degeneration types
fail at different stages of the solve, examined in turn below.

\textbf{Inverted meshes fail at assembly.} At the first inverted
element the Jacobian determinant becomes negative (for example
$\det J = -5.13 \times 10^{-4}$ at $\rho_d = 1\%$), and the
shape-function gradient routine raises a runtime error before any
linear system is assembled. To obtain a field for visual comparison
on the degenerate mesh (\cref{fig:degenerate}(a1)), we forced a
diagnostic solve on the $\rho_d = 25\%$ mesh with an
absolute-determinant Jacobian patch. CG did not converge after
$5{,}000$ iterations and the final relative residual was
$1.28 \times 10^6$, so the field is not a physical solution. Its
relative $L^2$ error is $99.9\%$, and its peak magnitude of
$2.12 \times 10^{-4}$~mm is only $0.39\%$ of the shared colour-scale
limit ($0.0540$~mm). The near-zero field in
\cref{fig:degenerate}(a2) reads as solver breakdown, not as a stiff
structural response.

\textbf{Sliver meshes fail one stage later, at the solve.} Because
$\det J$ stays positive at every Gauss point, the stiffness matrix
assembles, but the slivered elements (\cref{fig:degenerate}(b1))
raise the condition number until CG fails to converge within
$2{,}000$ iterations even at $\rho_d = 1\%$. A patched solve that
replaces $\det J$ with $|\det J|$ on the $\rho_d = 25\%$ mesh
returns the field in \cref{fig:degenerate}(b2). Its peak magnitude
($0.056$~mm) exceeds the reference maximum ($0.054$~mm) and
saturates the colour map, and its profile bears no resemblance to
the smooth axial gradient of a pressurised vessel. The relative
$L^2$ error is $63.2\%$, so bypassing the integrity check yields no
usable result either.

By contrast, FEVessel returns a smooth, physically consistent field on
every degenerate mesh (\cref{fig:degenerate}(a3)). For inverted
elements the mean relative $L^2$ error stays at the clean-mesh
baseline of $1.50\% \pm 0.62\%$ across all six ratios, the flat
curve in \cref{fig:degenerate}(c), and the error map shows no
localised defect at the inverted elements
(\cref{fig:degenerate}(a4)). The reason is structural and lies in
the Transolver attention. The physics-attention module first
aggregates all points into a small set of tokens by a weighted sum
(\cref{eq:slice}). That aggregation sums over the points, so it is
invariant to their ordering, and the self-attention then acts on
the tokens alone. The operator therefore reads only the per-point
coordinates and physical channels, never the element connectivity,
and it never evaluates an element Jacobian. Inverting an element
changes its connectivity but not its node coordinates, so the
operator input is unchanged and an inverted mesh is
indistinguishable from its valid parent. Removing the Jacobian
determinant from the forward graph eliminates this entire class of
failure by construction.

The inverted case is lossless because it leaves the node
coordinates untouched. Sliver degeneration instead perturbs those
coordinates, and therefore the FEVessel input. The mean error rises
monotonically with the sliver ratio, from $1.54\% \pm 0.62\%$ at
$\rho_d = 1\%$ to $3.66\% \pm 0.83\%$ at $\rho_d = 75\%$, about
$2.4\times$ the clean-mesh baseline (\cref{fig:degenerate}(d) and
\cref{tab:degenerate}). The prediction stays smooth
(\cref{fig:degenerate}(b3)) and the error concentrates where the
coordinates are most distorted (\cref{fig:degenerate}(b4)). Even
with $75\%$ of the elements ($23{,}652$ of $31{,}536$) severely
distorted, FEVessel stays below $4\%$, within the range usable for
early-stage design.

\begin{table}[!t]
  \centering
  \caption{Performance of FEM and FEVessel on degenerate meshes at
    varying degeneration ratios, averaged over $30$ test vessels.
    ``Fail'' indicates that the conventional FEM solver did not
    converge or could not assemble the stiffness matrix. Forced FEM
    errors at $\rho_d=25\%$ are diagnostic non-converged fields used
    in \cref{fig:degenerate}. FEVessel errors are the relative $L^2$ error of $\boldsymbol{u}$,
    reported as mean $\pm$ standard deviation.}
  \label{tab:degenerate}
  \setlength{\tabcolsep}{3pt}
  \begin{tabular*}{\linewidth}{@{\extracolsep{\fill}}llccc@{}}
    \toprule
    Type & Ratio & FEM status & FEM error & FEVessel error \\
    \midrule
    Inverted & $1\%$  & Fail (neg.\ Jacobian) & --- & $1.50\% \pm 0.62\%$ \\
    Inverted & $5\%$  & Fail (neg.\ Jacobian) & --- & $1.50\% \pm 0.62\%$ \\
    Inverted & $10\%$ & Fail (neg.\ Jacobian) & --- & $1.50\% \pm 0.62\%$ \\
    Inverted & $25\%$ & Fail (neg.\ Jacobian) & $99.9\%$ (forced) & $1.50\% \pm 0.62\%$ \\
    Inverted & $50\%$ & Fail (neg.\ Jacobian) & --- & $1.50\% \pm 0.62\%$ \\
    Inverted & $75\%$ & Fail (neg.\ Jacobian) & --- & $1.50\% \pm 0.62\%$ \\
    \midrule
    Sliver   & $1\%$  & Fail (CG diverged)    & --- & $1.54\% \pm 0.62\%$ \\
    Sliver   & $5\%$  & Fail (CG diverged)    & --- & $1.66\% \pm 0.61\%$ \\
    Sliver   & $10\%$ & Fail (CG diverged)    & --- & $1.81\% \pm 0.59\%$ \\
    Sliver   & $25\%$ & Fail (CG diverged)    & $63.2\%$ (forced) & $2.25\% \pm 0.55\%$ \\
    Sliver   & $50\%$ & Fail (CG diverged)    & --- & $2.96\% \pm 0.60\%$ \\
    Sliver   & $75\%$ & Fail (CG diverged)    & --- & $3.66\% \pm 0.83\%$ \\
    \bottomrule
  \end{tabular*}
\end{table}

Across the inverted and sliver cases together, FEVessel holds the
relative $L^2$ error between $1.50\%$ and $3.66\%$ over all twelve
configurations, while conventional FEM returns no solution on any
of them. To our knowledge this is the first systematic study of
mesh-independent solution on industrially relevant 3D pressure vessels with degenerate meshes.
It lets FEVessel act as a fallback solver whenever the meshing pipeline produces
a low-quality mesh, removing the manual remeshing that bottlenecks
design iteration on such geometries.

\section{Discussion}\label{sec:discussion}

The first experiment shows that a single FEVessel model covers
material, geometry, and boundary variation at once, with no
retraining for each new configuration. In early pressure vessel
design, each candidate is evaluated in milliseconds, so large
design-of-experiments sweeps become practical at an accuracy
sufficient for trade-off studies.

The second experiment shows that the FEVessel prediction, used
as the initial guess for the iterative solver, speeds up
convergence markedly without any loss of accuracy. The
convergence rate is still set by the chosen iterative algorithm,
and the prediction only changes the initial value. The final
solution therefore remains accurate in the finite element sense.

The third experiment shows that FEVessel's prediction accuracy
stays nearly constant for input densities above roughly $30\%$
of the training resolution. One trained model can therefore
serve a range of mesh resolutions, which lightens the meshing
burden in pre-processing.

The fourth experiment shows that FEVessel gives usable solutions
on degenerate meshes where conventional FEM fails entirely.
Label-free training needs no reference solution, and such meshes
cannot produce one. This is exactly where the value of the
framework lies. Mesh quality is then no longer a gate between the
CAD model and analysis, and an end-to-end automated pipeline
becomes possible.

\section{Conclusion}\label{sec:conclusion}

We have presented FEVessel, which extends the Pretrained Finite
Element Method from two dimensions to 3D pressure vessel
analysis, with Transolver as the backbone network. Focusing on
the capabilities
that matter most in engineering, we carried out four experiments.
These cover FEVessel's joint generalisation across material,
geometry, and boundary conditions, its warm-start acceleration,
cross-resolution prediction from coarse input, and solvability
on degenerate meshes. A single label-free model
reaches a $1.35\%$ relative displacement error and a $2.07\%$
strain error on out-of-sample vessels. On the commonly used AMG solver, it gives about a
$9.2\times$ end-to-end speedup while preserving engineering
accuracy. It further holds about a $3\%$ error at only $30\%$
of the training point density. It also gives usable solutions
on low-quality meshes where conventional FEM diverges entirely.

Several limitations remain. The present evaluation focuses on
static structural analysis, yet transient thermo-mechanical and
fatigue analyses are equally critical in pressure vessel design.
These also fit naturally with finite-difference temporal schemes
and form a worthwhile direction. FEVessel gives smooth
predictions on degenerate meshes, but its accuracy is still
bounded by the training distribution. Warm-start refinement on
degenerate meshes also deserves further study. We hope to combine
FEVessel with generative models for the inverse design of
pressure vessel topology. We also aim to deploy it as a real-time
analysis backbone for operational monitoring.

We believe that the combination of label-free training and mesh
independence will gradually change how pressure vessel analysis
is carried out. Mesh quality has long been the bottleneck between
the CAD model and numerical analysis. Point-cloud operators
instead let predictions sit directly on the raw geometry. Once
analysis is no longer limited by mesh quality, large numbers of
design candidates can be evaluated automatically and reliably.
This work brings that value of the Pretrained Finite Element Method to the 3D setting, and we hope it is one
step along this path.



\bmsection*{Author Contributions}

\textbf{Yipin Sun:} Methodology, Software, Validation, Investigation,
Data curation, Writing -- original draft, Visualization.
\textbf{Baiyang Zheng:} Software, Validation, Investigation,
Data curation.
\textbf{Yizheng Wang:} Conceptualization, Methodology, Supervision,
Writing -- review \& editing.
\textbf{Yuzhou Lin:} Methodology (pressure vessel design),
Investigation.
\textbf{Xiaoying Zhuang:} Conceptualization, Supervision, Funding
acquisition, Resources, Writing -- review \& editing.
\textbf{Timon Rabczuk:} Conceptualization, Supervision, Project
administration, Writing -- review \& editing.

\bmsection*{Acknowledgments}
The first author thanks Yizheng Wang and Timon Rabczuk for
their patient guidance throughout this work, and acknowledges
the financial support of the China National Scholarship.

\bmsection*{Conflicts of Interest}
The authors declare that they have no known competing financial
interests or personal relationships that could have appeared to
influence the work reported in this paper.

\bmsection*{Data Availability Statement}
The dataset generated in this work (approximately 13~GB) is publicly
available at
\url{https://1drv.ms/f/c/37600f10be213315/IgCW58XT-wF0Q4kdX0C5YuA7ASjJ69f9BGQe6CFsLwHtr0s}.
The source code will be released at
\url{https://github.com/YPG-8-P1Guan/pfem-fevessel} upon acceptance.

\bmsection*{ORCID}
\noindent\textit{Yipin Sun}~\orchid{https://orcid.org/0009-0003-7752-8321}~\url{https://orcid.org/0009-0003-7752-8321}\\
\textit{Baiyang Zheng}~\orchid{https://orcid.org/0000-0001-8616-6180}~\url{https://orcid.org/0000-0001-8616-6180}\\
\textit{Yizheng Wang}~\orchid{https://orcid.org/0000-0002-3899-7008}~\url{https://orcid.org/0000-0002-3899-7008}\\
\textit{Yuzhou Lin}~\orchid{https://orcid.org/0009-0006-5541-3875}~\url{https://orcid.org/0009-0006-5541-3875}\\
\textit{Xiaoying Zhuang}~\orchid{https://orcid.org/0000-0001-6562-2618}~\url{https://orcid.org/0000-0001-6562-2618}\\
\textit{Timon Rabczuk}~\orchid{https://orcid.org/0000-0002-7150-296X}~\url{https://orcid.org/0000-0002-7150-296X}

\bibliography{refs}

\begin{thebibliography}{32}
\providecommand{\natexlab}[1]{#1}
\providecommand{\url}[1]{\texttt{#1}}
\expandafter\ifx\csname urlstyle\endcsname\relax
  \providecommand{\doi}[1]{doi: #1}\else
  \providecommand{\doi}{doi: \begingroup \urlstyle{rm}\Url}\fi

\bibitem[Moss(2004)]{moss2004pressure}
Dennis~R Moss.
\newblock \emph{Pressure vessel design manual}.
\newblock Elsevier, 2004.

\bibitem[Peters et~al.(2023)Peters, Subramanian, and
  Sarzynski]{peters2023overview}
Daniel~T Peters, Kannan Subramanian, and Melanie Sarzynski.
\newblock Overview of revisions to the asme boiler and pressure vessel code
  section viii division 3 for the 2023 edition and near future.
\newblock In \emph{Pressure Vessels and Piping Conference}, volume 87479, page
  V004T05A006. American Society of Mechanical Engineers, 2023.

\bibitem[Alves et~al.(2025)Alves, Novais, da~Silva, and
  Junior]{alves2025probabilistic}
Mariana~Pimenta Alves, Henrique~Cordeiro Novais, Samuel da~Silva, and Carlos
  Alberto~Cimini Junior.
\newblock Probabilistic prediction of burst pressure from composite pressure
  vessels: Effect of material uncertainty and stacking sequence.
\newblock \emph{Results in Engineering}, page 106111, 2025.

\bibitem[Li et~al.(2025)Li, Lv, Zhang, He, and Zhang]{li2025reliability}
Wenbo Li, Hong Lv, Lijun Zhang, Pengfei He, and Cunman Zhang.
\newblock Reliability analysis and optimization design of hydrogen storage
  composite pressure vessel with hybrid random-fuzzy uncertainties.
\newblock \emph{Journal of Reinforced Plastics and Composites}, 44\penalty0
  (7-8):\penalty0 358--374, 2025.

\bibitem[Johnson et~al.(2023)Johnson, Zhu, Sindelar, and
  Wiersma]{johnson2023parametric}
William~R Johnson, Xian-Kui Zhu, Robert Sindelar, and Bruce Wiersma.
\newblock A parametric finite element study for determining burst strength of
  thin and thick-walled pressure vessels.
\newblock \emph{International Journal of Pressure Vessels and Piping},
  204:\penalty0 104968, 2023.

\bibitem[Zhang et~al.(2022)Zhang, Yu, and Liu]{zhang2022surrogate}
Long Zhang, Chengguo Yu, and Bingbin Liu.
\newblock Surrogate-based structural optimization design of large-scale
  rectangular pressure vessel using radial point interpolation method.
\newblock \emph{International Journal of Pressure Vessels and Piping},
  197:\penalty0 104638, 2022.

\bibitem[Kang et~al.(2025)Kang, Kang, Kim, and Lee]{kang2025condition}
Dayoung Kang, Shinseong Kang, Bongseok Kim, and Kyunghoon Lee.
\newblock Condition-based fatigue life monitoring of a high-pressure hydrogen
  storage vessel using a reduced basis digital twin.
\newblock \emph{Engineering Structures}, 336:\penalty0 120196, 2025.

\bibitem[Wang et~al.(2024)Wang, Sun, Rabczuk, and Liu]{wang2024dcem}
Yizheng Wang, Jia Sun, Timon Rabczuk, and Yinghua Liu.
\newblock Dcem: A deep complementary energy method for linear elasticity.
\newblock \emph{International Journal for Numerical Methods in Engineering},
  125\penalty0 (24):\penalty0 e7585, 2024.

\bibitem[Zienkiewicz et~al.(1977)Zienkiewicz, Taylor, Nithiarasu, and
  Zhu]{zienkiewicz1977finite}
Olgierd~Cecil Zienkiewicz, Robert~Leroy Taylor, Perumal Nithiarasu, and JZ~Zhu.
\newblock \emph{The finite element method}, volume~3.
\newblock Elsevier, 1977.

\bibitem[Hughes(2003)]{hughes2003finite}
Thomas~JR Hughes.
\newblock \emph{The finite element method: linear static and dynamic finite
  element analysis}.
\newblock Courier Corporation, 2003.

\bibitem[Bozkurt et~al.(2021)Bozkurt, Nash, and Uzzaman]{bozkurt2021comparison}
Murat Bozkurt, David Nash, and Asraf Uzzaman.
\newblock A comparison of stress analysis and limit analysis approaches for
  single and multiple nozzle combinations in cylindrical pressure vessels.
\newblock \emph{International Journal of Pressure Vessels and Piping},
  194:\penalty0 104563, 2021.

\bibitem[Wang et~al.(2026)Wang, Hao, Eshaghi, Anitescu, Zhuang, Rabczuk, and
  Liu]{wang2026pretrain}
Yizheng Wang, Zhongkai Hao, Mohammad~Sadegh Eshaghi, Cosmin Anitescu, Xiaoying
  Zhuang, Timon Rabczuk, and Yinghua Liu.
\newblock Pretrain finite element method: A pretraining and warm-start
  framework for pdes via physics-informed neural operators.
\newblock \emph{Journal of the Mechanics and Physics of Solids}, 214:\penalty0
  106682, 2026.

\bibitem[Kudela and Matousek(2022)]{kudela2022recent}
Jakub Kudela and Radomil Matousek.
\newblock Recent advances and applications of surrogate models for finite
  element method computations: a review: J. kudela, r. matousek.
\newblock \emph{Soft Computing}, 26\penalty0 (24):\penalty0 13709--13733, 2022.

\bibitem[Sun et~al.(2025)Sun, Wen, Li, Cao, and Fei]{sun2025novel}
Yi-Pin Sun, Jiong-Ran Wen, Jian Li, Ai-Fang Cao, and Cheng-Wei Fei.
\newblock Novel integrated model approach for high cycle fatigue life and
  reliability assessment of helicopter flange structures.
\newblock \emph{Aerospace}, 12\penalty0 (2):\penalty0 78, 2025.

\bibitem[Flanagan et~al.(2025)Flanagan, O’Connor, Erfanian, Music, Brambley,
  and O’Kiely]{flanagan2025through}
Francis Flanagan, Alison~N O’Connor, Mozhdeh Erfanian, Omer Music, Edward~J
  Brambley, and Doireann O’Kiely.
\newblock Through-thickness resolution, stress oscillations and residual stress
  in cold rolling.
\newblock \emph{European Journal of Mechanics-A/Solids}, 114:\penalty0 105761,
  2025.

\bibitem[Satheesh et~al.(2024)Satheesh, Schmidt, Wall, and
  Meier]{satheesh2024novel}
Abhiroop Satheesh, Christoph~P Schmidt, Wolfgang~A Wall, and Christoph Meier.
\newblock A novel mesh regularization approach based on finite element
  distortion potentials: Application to material expansion processes with
  extreme volume change.
\newblock \emph{Computer Methods in Applied Mechanics and Engineering},
  432:\penalty0 117444, 2024.

\bibitem[van Huyssteen and Reddy(2022)]{van2022incorporation}
Daniel van Huyssteen and BD~Reddy.
\newblock The incorporation of mesh quality in the stabilization of virtual
  element methods for nonlinear elasticity.
\newblock \emph{Computer Methods in Applied Mechanics and Engineering},
  392:\penalty0 114720, 2022.

\bibitem[Wei et~al.(2024)Wei, Lei, Pei, and Dong]{wei2024structural}
Zhigang Wei, Linsen Lei, Xianjun Pei, and Pingsha Dong.
\newblock The structural strain method for fatigue evaluation of welded
  components: Analytical treatment of reversed plasticity.
\newblock \emph{International Journal of Pressure Vessels and Piping},
  210:\penalty0 105249, 2024.

\bibitem[Sukumar and Tupek(2022)]{sukumar2022virtual}
N~Sukumar and Michael~R Tupek.
\newblock Virtual elements on agglomerated finite elements to increase the
  critical time step in elastodynamic simulations.
\newblock \emph{International Journal for Numerical Methods in Engineering},
  123\penalty0 (19):\penalty0 4702--4725, 2022.

\bibitem[Prabhune and Suresh(2024)]{prabhune2024isoparametric}
Bhagyashree Prabhune and Krishnan Suresh.
\newblock An isoparametric tangled finite element method for handling
  higher-order elements with negative jacobian.
\newblock \emph{Computational Mechanics}, 73\penalty0 (1):\penalty0 159--176,
  2024.

\bibitem[Conley et~al.(2016)Conley, Delaney, and Jiao]{conley2016overcoming}
Rebecca Conley, Tristan~J Delaney, and Xiangmin Jiao.
\newblock Overcoming element quality dependence of finite elements with
  adaptive extended stencil fem (aes-fem).
\newblock \emph{International Journal for Numerical Methods in Engineering},
  108\penalty0 (9):\penalty0 1054--1085, 2016.

\bibitem[Karniadakis et~al.(2021)Karniadakis, Kevrekidis, Lu, Perdikaris, Wang,
  and Yang]{karniadakis2021physics}
George~Em Karniadakis, Ioannis~G Kevrekidis, Lu~Lu, Paris Perdikaris, Sifan
  Wang, and Liu Yang.
\newblock Physics-informed machine learning.
\newblock \emph{Nature Reviews Physics}, 3\penalty0 (6):\penalty0 422--440,
  2021.

\bibitem[Cuomo et~al.(2022)Cuomo, Di~Cola, Giampaolo, Rozza, Raissi, and
  Piccialli]{cuomo2022scientific}
Salvatore Cuomo, Vincenzo~Schiano Di~Cola, Fabio Giampaolo, Gianluigi Rozza,
  Maziar Raissi, and Francesco Piccialli.
\newblock Scientific machine learning through physics--informed neural
  networks: Where we are and what’s next.
\newblock \emph{Journal of Scientific Computing}, 92\penalty0 (3):\penalty0 88,
  2022.

\bibitem[Raissi et~al.(2019)Raissi, Perdikaris, and
  Karniadakis]{raissi2019physics}
Maziar Raissi, Paris Perdikaris, and George~E Karniadakis.
\newblock Physics-informed neural networks: A deep learning framework for
  solving forward and inverse problems involving nonlinear partial differential
  equations.
\newblock \emph{Journal of Computational physics}, 378:\penalty0 686--707,
  2019.

\bibitem[Lu et~al.(2021)Lu, Jin, Pang, Zhang, and Karniadakis]{lu2021learning}
Lu~Lu, Pengzhan Jin, Guofei Pang, Zhongqiang Zhang, and George~Em Karniadakis.
\newblock Learning nonlinear operators via deeponet based on the universal
  approximation theorem of operators.
\newblock \emph{Nature machine intelligence}, 3\penalty0 (3):\penalty0
  218--229, 2021.

\bibitem[Li et~al.(2021)Li, Kovachki, Azizzadenesheli, Liu, Bhattacharya,
  Stuart, and Anandkumar]{li2021fourier}
Zongyi Li, Nikola~Borislavov Kovachki, Kamyar Azizzadenesheli, Burigede Liu,
  Kaushik Bhattacharya, Andrew Stuart, and Anima Anandkumar.
\newblock Fourier neural operator for parametric partial differential
  equations.
\newblock In \emph{International Conference on Learning Representations}, 2021.

\bibitem[Li et~al.(2024)Li, Zheng, Kovachki, Jin, Chen, Liu, Azizzadenesheli,
  and Anandkumar]{li2024physics}
Zongyi Li, Hongkai Zheng, Nikola Kovachki, David Jin, Haoxuan Chen, Burigede
  Liu, Kamyar Azizzadenesheli, and Anima Anandkumar.
\newblock Physics-informed neural operator for learning partial differential
  equations.
\newblock \emph{ACM/IMS Journal of Data Science}, 1\penalty0 (3):\penalty0
  1--27, 2024.

\bibitem[Eshaghi et~al.(2025)Eshaghi, Anitescu, Thombre, Wang, Zhuang, and
  Rabczuk]{eshaghi2025variational}
Mohammad~Sadegh Eshaghi, Cosmin Anitescu, Manish Thombre, Yizheng Wang,
  Xiaoying Zhuang, and Timon Rabczuk.
\newblock Variational physics-informed neural operator (vino) for solving
  partial differential equations.
\newblock \emph{Computer Methods in Applied Mechanics and Engineering},
  437:\penalty0 117785, 2025.

\bibitem[Wu et~al.(2024)Wu, Luo, Wang, Wang, and Long]{Wu2024TransolverAF}
Haixu Wu, Huakun Luo, Haowen Wang, Jianmin Wang, and Mingsheng Long.
\newblock Transolver: A fast transformer solver for pdes on general geometries.
\newblock In \emph{International Conference on Machine Learning}, 2024.

\bibitem[Timoshenko and Goodier(1970)]{timoshenko1970theory}
Stephen~P Timoshenko and James~N Goodier.
\newblock \emph{Theory of elasticity}.
\newblock McGraw-Hill, 3 edition, 1970.

\bibitem[Saad(2003)]{saad2003iterative}
Yousef Saad.
\newblock \emph{Iterative methods for sparse linear systems}.
\newblock SIAM, 2003.

\bibitem[Van{\v{e}}k et~al.(1996)Van{\v{e}}k, Mandel, and
  Brezina]{vanek1996algebraic}
Petr Van{\v{e}}k, Jan Mandel, and Marian Brezina.
\newblock Algebraic multigrid by smoothed aggregation for second and fourth
  order elliptic problems.
\newblock \emph{Computing}, 56\penalty0 (3):\penalty0 179--196, 1996.

\end{thebibliography}

\appendix
\section{Qualitative displacement montage across the
M/G/B parameter space}\label{app:mgb-montage}
\begin{figure*}[p]
  \centering
  \includegraphics[width=\linewidth]{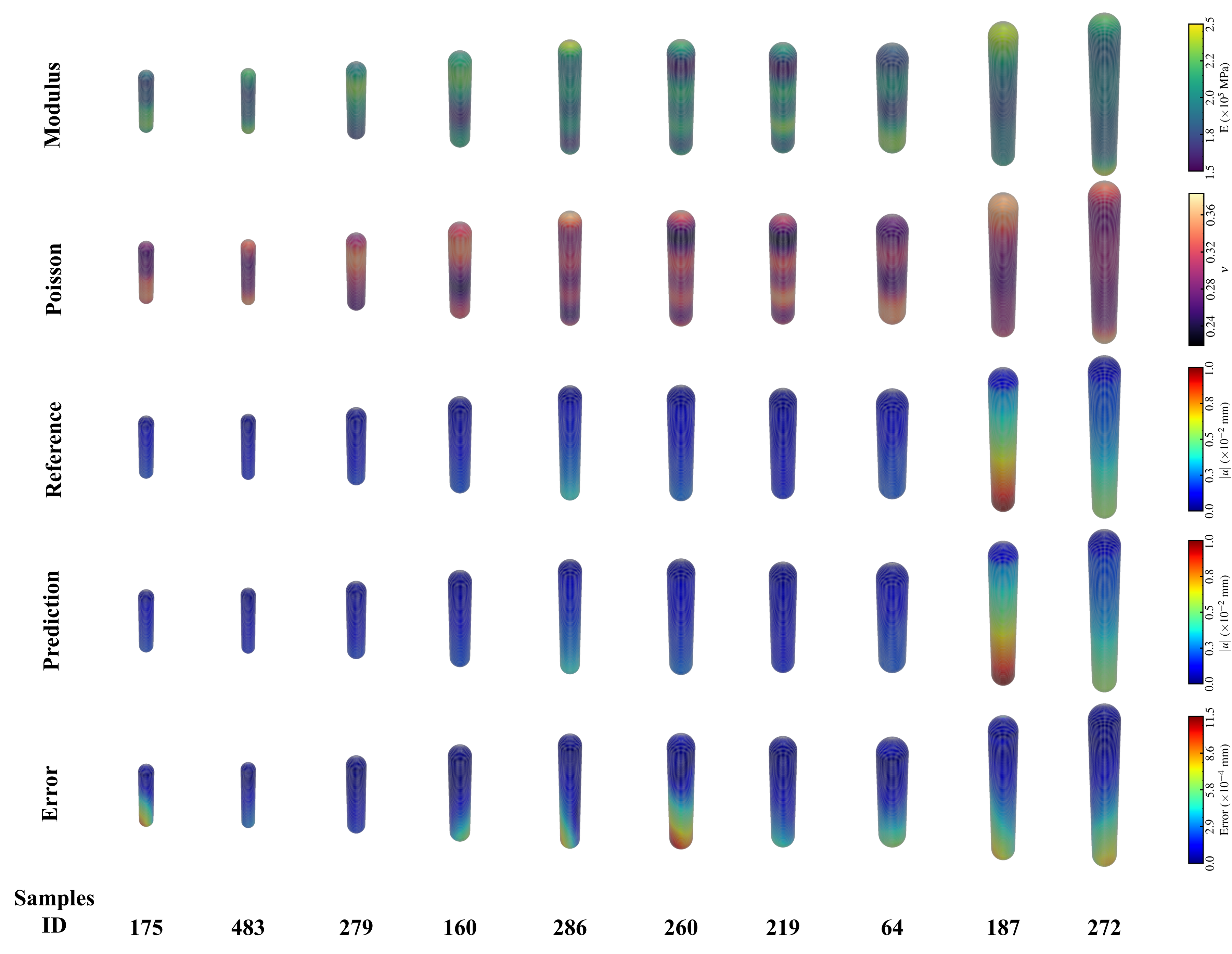}
  \caption{Predicted displacements and pointwise vector-error
    contours for 10 test cases across the extremes of the
    M/G/B space. Rows, top to bottom:
    modulus $E$, Poisson's ratio $\nu$, FEM reference
    $|\boldsymbol{u}_{\text{ref}}|$, FEVessel prediction
    $|\boldsymbol{u}_{\text{pred}}|$, and vector error
    $\|\boldsymbol{u}_{\text{pred}}-\boldsymbol{u}_{\text{ref}}\|_2$;
    sample IDs at the bottom.}
  \label{fig:mgb-contour}
\end{figure*}

\Cref{fig:mgb-contour} extends the representative comparison of
\cref{subsec:mgb} to ten test vessels spanning the extremes of the
M/G/B parameter space, confirming that the visual agreement
between FEVessel and the FEM reference holds across the full range of
materials, geometries, and pressure boundaries.

Across these ten cases, the FEVessel predictions
are visually indistinguishable from the FEM reference for the
full range of vessel sizes, aspect ratios, and pressure boundary
configurations: tall slender vessels (e.g.\ Sample~IDs 175 and 219),
short stubby vessels (Sample~IDs 286 and 483), and partial-pressure
cases where the loaded surface covers only a fraction of the inner
cavity (Sample~IDs 64 and 272) are all captured with the same
fidelity. The largest local errors, on the order of
$8\times10^{-4}$~mm, concentrate near the hemispherical end-cap
transitions, where displacement gradients are largest, consistent
with the expectation that stress-concentration regions are the
hardest for any operator-learning model to approximate.

\section{Convergence acceleration of the warm
start}\label{app:warmstart-theory}

\Cref{subsec:warmstart} reports that the FEVessel warm start lowers the
iteration count of every solver tested. This appendix derives that
behaviour for the linear pressure vessel system. The argument follows
the convergence analysis in Appendix~C of Wang
et al.~\cite{wang2026pretrain}, specialised here to the symmetric
positive definite (SPD) stiffness system and extended with a
wall-clock corollary that accounts for the solver setup cost.

The assembled system $\mathbf{K}\mathbf{U}=\mathbf{F}$ of
\cref{eq:warmstart} is SPD, because the Dirichlet constraint on
$\Gamma_u$ removes the rigid-body modes (\cref{subsec:problem}), so a
unique solution $\mathbf{U}^{\star}=\mathbf{K}^{-1}\mathbf{F}$ exists.
For the iterate $\mathbf{U}^{(k)}$ we define the error and the
residual,
\begin{equation}\label{eq:app-err-res}
  \mathbf{e}^{(k)}=\mathbf{U}^{(k)}-\mathbf{U}^{\star},
  \qquad
  \mathbf{r}^{(k)}=\mathbf{K}\mathbf{U}^{(k)}-\mathbf{F}
  =\mathbf{K}\,\mathbf{e}^{(k)},
\end{equation}
with the residual sign matching \cref{eq:residual}. We use the energy
norm $\lVert\mathbf{v}\rVert_{\mathbf{K}}
=(\mathbf{v}^{\top}\mathbf{K}\mathbf{v})^{1/2}$ and write
$\kappa=\lambda_{\max}/\lambda_{\min}$ for the spectral condition
number of $\mathbf{K}$, with $\lambda_{\min}$ and $\lambda_{\max}$ the
extreme eigenvalues. The residual norm and the energy-norm error are
linked by
\begin{equation}\label{eq:app-norm-link}
  \lambda_{\min}\lVert\mathbf{e}^{(k)}\rVert
  \le\lVert\mathbf{r}^{(k)}\rVert
  \le\lambda_{\max}\lVert\mathbf{e}^{(k)}\rVert ,
\end{equation}
so geometric decay of the error implies geometric decay of the
residual monitored by \cref{eq:residual}.

We first take the algorithm $\boldsymbol{\phi}$ of
\cref{eq:warmstart} to be a stationary iteration built on a
preconditioner $\mathbf{B}\approx\mathbf{K}^{-1}$,
\begin{equation}\label{eq:app-stationary}
  \mathbf{U}^{(k+1)}=\mathbf{U}^{(k)}-\mathbf{B}\,\mathbf{r}^{(k)} ,
\end{equation}
of which a single algebraic multigrid V-cycle is an instance with
$\mathbf{B}$ the multigrid operator. Subtracting $\mathbf{U}^{\star}$
from \cref{eq:app-stationary} and using
$\mathbf{F}-\mathbf{K}\mathbf{U}^{\star}=\mathbf{0}$, the error obeys
\begin{equation}\label{eq:app-error-recursion}
  \mathbf{e}^{(k+1)}=(\mathbf{I}-\mathbf{B}\mathbf{K})\,\mathbf{e}^{(k)}
  =\mathbf{G}\,\mathbf{e}^{(k)},
  \qquad
  \mathbf{e}^{(k)}=\mathbf{G}^{k}\mathbf{e}^{(0)} ,
\end{equation}
which is the recursion of Wang et al.~\cite{wang2026pretrain} written
with the iteration matrix $\mathbf{G}=\mathbf{I}-\mathbf{B}\mathbf{K}$.
We denote the energy-norm contraction factor by
$\rho=\lVert\mathbf{G}\rVert_{\mathbf{K}}$. For a convergent method
$\rho<1$, and for smoothed-aggregation AMG supplied with the
rigid-body near-null-space modes $\rho$ is bounded below one
independently of the node count under standard
assumptions~\cite{vanek1996algebraic}. Taking the energy norm of
\cref{eq:app-error-recursion} gives
\begin{equation}\label{eq:app-contraction}
  \lVert\mathbf{e}^{(k)}\rVert_{\mathbf{K}}
  \le\rho^{k}\,\lVert\mathbf{e}^{(0)}\rVert_{\mathbf{K}} .
\end{equation}
The contraction factor $\rho$ depends only on $\mathbf{K}$ and the
solver, not on the initial guess. The initial guess enters
\cref{eq:app-contraction} only through the prefactor
$\lVert\mathbf{e}^{(0)}\rVert_{\mathbf{K}}$, so a better initial guess
lowers the whole convergence curve by a constant factor and leaves
its slope unchanged.

The stopping rule of \cref{eq:residual} acts on the relative residual
$\eta^{(k)}=\lVert\mathbf{r}^{(k)}\rVert/\lVert\mathbf{F}\rVert$. The
residual obeys the same recursion as the error and contracts at the
same factor $\rho$, so $\eta^{(k)}\le\rho^{k}\eta^{(0)}$ holds up to a
constant from \cref{eq:app-norm-link} that does not depend on the
tolerance. Convergence is reached once $\rho^{k}\eta^{(0)}\le
\mathrm{tol}$, which gives the iteration-count estimate
\begin{equation}\label{eq:app-itercount}
  n\!\left(\mathbf{U}^{(0)}\right)
  =\left\lceil
    \frac{\ln\!\left(\eta^{(0)}/\mathrm{tol}\right)}{\ln(1/\rho)}
  \right\rceil .
\end{equation}
This estimate specialises to the relative residual the result of Wang
et al.~\cite{wang2026pretrain}, who show that a smaller initial error
reduces the upper bound on the iteration count. The two
initialisations of \cref{subsec:warmstart} give $\eta^{(0)}=1$ for the
zero guess $\mathbf{U}^{(0)}=\mathbf{0}$, since
$\mathbf{r}^{(0)}=-\mathbf{F}$, and $\eta^{(0)}=\eta_{0}$ for the warm
start $\mathbf{U}^{(0)}=\mathbf{U}^{\mathrm{NO}}$, with
$\eta_{0}\approx0.021$ measured over the test set. Subtracting the two
counts from \cref{eq:app-itercount}, where the tolerance-independent
constant cancels, the warm start saves
\begin{equation}\label{eq:app-deltan}
  \Delta n
  =n_{\text{zero}}-n_{\text{warm}}
  =\frac{\ln(1/\eta_{0})}{\ln(1/\rho)}
  =\frac{-\log_{10}\eta_{0}}{\log_{10}(1/\rho)} .
\end{equation}
Three properties of \cref{eq:app-deltan} match the measurements of
\cref{subsec:warmstart}. First, $\eta_{0}<1$ and $\rho<1$ give
$\Delta n>0$, and $\Delta n$ grows with the residual reduction factor
$1/\eta_{0}$, so a closer warm start saves more iterations. Second,
$\Delta n$ carries no tolerance, so the warm start removes a fixed
block of iterations set by the initial-residual reduction
$-\log_{10}\eta_{0}\approx1.68$ decades divided by the per-iteration
decade rate. This fixed reduction of about $1.68$ decades is the value
reported in \cref{subsec:warmstart} and drawn in
\cref{fig:warmstart-conv}(d). Third, the share of the total residual
budget supplied by the warm start,
\begin{equation}\label{eq:app-share}
  s=\frac{\Delta n}{n_{\text{zero}}}
  =\frac{\log_{10}(1/\eta_{0})}{\log_{10}(1/\mathrm{tol})} ,
\end{equation}
falls as the tolerance tightens, from about $84\%$ at
$\mathrm{tol}=10^{-2}$ to $56\%$ at $10^{-3}$ and $42\%$ at
$10^{-4}$. These shares reproduce the decomposition of
\cref{fig:warmstart-conv}(d) and recover the conclusion of Wang
et al.~\cite{wang2026pretrain} that the benefit of a good initial
guess diminishes at very small tolerances.

When $\boldsymbol{\phi}$ is the conjugate gradient method or its
AMG-preconditioned form, the iterate minimises the energy-norm error
over the Krylov subspace, giving the standard
bound~\cite{saad2003iterative}
\begin{equation}\label{eq:app-cg}
  \lVert\mathbf{e}^{(k)}\rVert_{\mathbf{K}}
  \le2\left(\frac{\sqrt{\kappa}-1}{\sqrt{\kappa}+1}\right)^{k}
  \lVert\mathbf{e}^{(0)}\rVert_{\mathbf{K}} ,
\end{equation}
with $\kappa=\kappa(\mathbf{K})$ for CG and
$\kappa=\kappa(\mathbf{B}\mathbf{K})$ for the preconditioned form. The
prefactor structure matches \cref{eq:app-contraction}, the rate
depending only on $\kappa$ and the initial guess entering only through
$\lVert\mathbf{e}^{(0)}\rVert_{\mathbf{K}}$. Repeating the argument of
\cref{eq:app-itercount,eq:app-deltan} with $1/\rho$ replaced by
$(\sqrt{\kappa}+1)/(\sqrt{\kappa}-1)$ yields the same
tolerance-independent saving, the constant $2$ and the residual
conversion of \cref{eq:app-norm-link} adding only a bounded offset.
The warm start therefore accelerates CG and AMG-PCG by the same
mechanism. A single contraction factor makes the saving exactly
fixed, whereas the early transient of a real solver contracts faster
than the asymptotic rate, so the measured iteration speedups at the
loosest tolerance in \cref{tab:warmstart} exceed this estimate while
following the same decreasing trend.

\subsection{End-to-end speedup with fixed setup
cost}\label{app:warmstart-walltime}

The saving of \cref{eq:app-deltan} counts iterations, whereas
\cref{subsec:warmstart} reports wall-clock speedups, and the two
differ because every AMG-based solve carries a fixed cost that the
warm start cannot remove. The per-solve wall-clock time splits into
\begin{equation}\label{eq:app-walltime}
  T\!\left(\mathbf{U}^{(0)}\right)
  =T_{\text{setup}}+n\!\left(\mathbf{U}^{(0)}\right)T_{\text{iter}} ,
\end{equation}
where $T_{\text{iter}}$ is the cost of one iteration and
$T_{\text{setup}}$ is the one-time setup cost. For AMG and
AMG-PCG, $T_{\text{setup}}$ builds the multigrid hierarchy on the
current stiffness matrix, and
the warm start adds one FEVessel forward pass to it. The defining
feature of $T_{\text{setup}}$ is that it does not depend on the
initial guess.
The zero guess and the warm start build the same hierarchy on the same
matrix, so the warm start lowers only the second term of
\cref{eq:app-walltime}. Writing $\gamma=T_{\text{setup}}/T_{\text{iter}}$
for the setup cost in iteration units and using
$\Delta n=n_{\text{zero}}-n_{\text{warm}}$, the wall-clock speedup is
\begin{equation}\label{eq:app-speedup}
  S(\mathrm{tol})
  =\frac{T_{\text{zero}}}{T_{\text{warm}}}
  =\frac{\gamma+n_{\text{zero}}}{\gamma+n_{\text{warm}}}
  =1+\frac{\Delta n}{\gamma+n_{\text{warm}}} .
\end{equation}
\Cref{eq:app-speedup} is an exact identity once the measured iteration
counts are inserted, and it separates the three solver behaviours of
\cref{subsec:warmstart}.

For unpreconditioned CG there is no hierarchy to build, yet the warm
start still adds one FEVessel forward pass and skips the high-residual
early iterations, so the per-iteration cost is not constant. Its
wall-clock speedup therefore follows the iteration saving closely
without equalling it, falling from $71.9\times$ at $10^{-2}$ to
$1.02\times$ at $10^{-4}$ against iteration speedups of $88.0\times$
and $1.05\times$.

For AMG the setup is not empty, $\gamma>0$, and \cref{eq:app-speedup}
throttles the speedup at both ends of the tolerance range. At the
loosest tolerance the warm start converges in barely one cycle, so
$n_{\text{warm}}$ is close to one and the denominator is dominated by
$\gamma$. The speedup is then capped near $1+\Delta n/\gamma$ by the
setup cost, well below the iteration speedup, which is why
\cref{tab:warmstart} records a $7.7\times$ wall-clock gain at
$10^{-2}$ against a $58.4\times$ iteration gain. At the tightest
tolerance $n_{\text{warm}}$ grows large, the denominator is dominated
by $n_{\text{warm}}$, and the speedup decays towards one, reaching
$2.7\times$ at $10^{-4}$. The maximum sits between these limits, at the
engineering tolerance $10^{-3}$, where the saving $\Delta n$ is still
large relative to $\gamma$ and the warm-start count has not yet
inflated the denominator. The end-to-end AMG speedup therefore peaks
at the tolerance prescribed by structural design practice rather than
at the loosest tolerance.

For AMG-PCG the setup cost is paid under both initialisations, as in
AMG, but the near-null-space preconditioner already keeps the
zero-start count small, so the saving $\Delta n$ is modest. The
speedup stays low and nearly flat across the three tolerances, between
$1.51\times$ and $1.55\times$ in \cref{tab:warmstart}, consistent with
a small $\Delta n$ over a denominator dominated by $\gamma$.

A single contraction factor would make $\Delta n$ constant and
\cref{eq:app-speedup} monotone in the tolerance. The non-monotonic AMG
peak instead reflects the measured counts, in which the early
transient lets the zero-start count grow faster than the warm-start
count between $10^{-2}$ and $10^{-3}$ and so enlarges $\Delta n$ there,
while the setup term $\gamma$ caps the loosest tolerance.
\Cref{eq:app-speedup} with the measured counts captures both effects.

\end{document}